\documentclass[12pt]{amsart}
\usepackage{pifont}
\usepackage{amsfonts}
\usepackage{amsmath}
\usepackage{amssymb}
\usepackage{amsmath,bm}
\usepackage{cite}
\usepackage{comment}
\usepackage{color}

\newtheorem{theorem}{Theorem}[section]

\newtheorem{lemma}[theorem]{Lemma}
\newtheorem{proposition}[theorem]{Proposition}

\theoremstyle{definition}
\newtheorem{definition}[theorem]{Definition}
\newtheorem{remark}{Remark}[section]

\newcommand{\R}{\mathbb{R}}

\begin{document}

\title[Generalized principal eigenvalue for nonlocal operators]{Generalized principal eigenvalues for nonlocal operators with drift}

\author[tong]{Hailin Tong}
\address{Hailin Tong,  School of Mathematical Sciences, Anhui University, Hefei 230601, Anhui, China}
\email{A24201022@stu.ahu.edu.cn}

\author[Zhou] {Tao Zhou$^*$}
\thanks{$^*$ Corresponding author: Tao Zhou}
\address{Tao Zhou, Center for Pure Mathematics, School of Mathematical Sciences, Anhui University, Hefei 230601, Anhui, China}
\email{tzhou910@ustc.edu.cn}

%\thanks{T. Zhou is supported by NSFC grant (12471150).}

\begin{abstract}
     Nonlocal diffusion operators with a drift term naturally arise in models such as moving coordinate systems induced by climate change. 
     This paper systematically studies several concepts of generalized principal eigenvalues (GPEs) for such operators. 
     We allow the dispersal kernel to be heterogeneous and non-symmetric, and remove the standard global upper bound on the kernel (thereby covering the case of heavy-tailed dispersal kernels).
     To overcome the analytical difficulties brought by this more general setting, we establish a new Harnack-type inequality for non-compactly supported kernels, and apply  homogenization and cut-off approximation techniques to derive two-sided bounds for the domain approximation limit of GPEs. 
     Furthermore, under the almost periodic setting, we prove that the GPEs defined through different admissible test function spaces (e.g., \(C^1\) spaces with bounded derivatives and uniformly continuous function spaces in the viscosity sense) are completely coincident. 
     This characterization reveals a fundamental unification: the limit result of principal eigenvalues under domain approximation by Coville and Hamel  coincides exactly with the explicit variational formula  by Du, Fang, and Ni, 
     thereby bridging two previously independent theoretical results.
\end{abstract}


\maketitle
\section{Introduction}

The theory of principal eigenvalues serves as a fundamental bridge between the spectral analysis of linear operators and nonlinear reaction-diffusion equations. For uniformly elliptic operators defined on bounded domains, the classical Krein--Rutman theorem guarantees the existence and simplicity of the principal eigenpair. The sign of this eigenvalue directly determines the existence of positive solutions to semilinear problems and the stability of the trivial steady state, thereby playing an indispensable role in persistence criteria, the construction of traveling wave solutions, and the characterization of maximum principles.

However, the classical local diffusion model 
$$u_t=\Delta u+f(u)$$
 assumes that individuals only move between adjacent spatial locations, and thus cannot describe long-range dispersal phenomena such as wind-borne seed dispersal or long-distance animal migration. 
 To address this limitation, researchers introduce a nonnegative kernel function \(J(x-y)\) to characterize the probability density of an individual migrating from \(y\) to \(x\), replacing the local Laplacian with the integral operator 
 $$\int J(x-y)u(y)\,dy-u(x).$$ 
 When the kernel \(J\) shrinks to the origin at an appropriate scale, the nonlocal model converges to the classical local diffusion equation, and can therefore be regarded as a natural generalization of the latter (see, for example, \cite{Coville2005Propagation}).

The introduction of a drift term into nonlocal diffusion models has clear applied motivations.
In studying the impact of climate change on population persistence, the environmentally favorable region typically shifts spatially at a constant speed \(c\), and the nonlinear term takes the form \(f(t,x,u)=\tilde f(x-ct,u)\). 
Upon transforming to the moving coordinate \(\xi=x-ct\), an explicit advective term \(c\,u_\xi\) naturally appears in the equation. 
Moreover, an asymmetric dispersal kernel \(J(x)\neq J(-x)\) itself implies directed transport, leading to unequal leftward and rightward spreading speeds. These facts naturally motivate the spectral study of nonlocal operators with drift, namely
\begin{equation}\label{eq:the-operator}
	\mathcal{L}^{\Omega}\phi(x) = \int_{\Omega} \phi(y) J(x,y) \, dy + a(x)\phi(x) + q(x)\phi'(x),
\end{equation}
and in particular the analysis of their principal eigenvalues.

Unlike local elliptic operators, nonlocal integral operators are generally not compact: when the potential \(a\) is nonconstant, neither the operator nor its resolvent is compact, so the classical Krein--Rutman theory cannot be applied directly; the operator may even possess no eigenvalues at all in standard function spaces (\cite[Section 6]{ShenZhang2010}).
 The introduction of an advective term further breaks symmetry and self-adjointness, rendering the spectral structure even more complex. 
 Against this background, the notion of a generalized principal eigenvalue becomes an indispensable theoretical tool. 
 The idea behind this concept traces back to the Collatz--Wielandt formula for positive matrices in the Perron--Frobenius theory, and was later extended to elliptic operators by Donsker and Varadhan \cite{DonskerVaradhan1975}, who defined
$$
\lambda(\Omega)=\inf\{\lambda:\exists\phi>0 \text{ in  some appropriate admissible set, }\ L\phi\le\lambda\phi \text{ in }\Omega\},
$$
where $L$ is a second order elliptic operator.
Berestycki, Nirenberg, and Varadhan \cite{BNV} subsequently proposed the complementary definition
$$
\tilde\lambda(\Omega)=\inf\{\lambda:\exists\phi>0 \text{ in  some appropriate admissible set, }\ L\phi\ge\lambda\phi \text{ in }\Omega\},
$$
and proved that the two notions coincide on smooth bounded domains, see also \cite{BR2015CPAM} for general domain.
 Coville \cite{coville2010jde,Coville2013}  extended this framework to nonlocal operators and systematically developed the theory of generalized principal eigenvalues for nonlocal operators.
In the purely nonlocal setting without drift, i.e. $q=0$ in \eqref{eq:the-operator}, Coville provided a sufficient condition for the existence of a continuous positive eigenfunction on bounded domains:
$$
\frac{1}{a_\infty-a(x)}\in L^1(\Omega).
$$
When this condition fails, the principal eigenvector can only be constructed in the space of measures, and is singular with respect to the Lebesgue measure. 
This finding demonstrates that, even on bounded domains, the principal eigenvalue of a nonlocal operator need not correspond to an ordinary function. 
Subsequently, Berestycki, Coville, and Vo \cite{BCV}, among others, further established Collatz--Wielandt-type characterizations for nonlocal operators, refining the theoretical foundation of generalized principal eigenvalues.

Once a drift term is incorporated into nonlocal operators, the spectral theory takes on an even richer landscape. 
Coville, Li and Wang \cite{LCW} were the first to study the one-dimensional case with a drift of fixed sign, proving the existence of a principal eigenpair on bounded domains with boundary conditions that must be consistent with the direction of advection (zero boundary value at the outflow end), and establishing a Collatz--Wielandt-type equivalent characterization. 
For unbounded domains, Coville and Hamel \cite{Coville2020} proved the crucial domain approximation property
$$
\lim_{n\to\infty}\lambda(\Omega_n)=\lambda(\Omega),
$$
where $\Omega_n\nearrow\Omega$. This reduces the problem on the unbounded domain $\Omega$ to a limit of problems on bounded truncated domains $\Omega_n$.
On the other hand, Du, Fang, and Ni \cite{DuFangNi2025} obtained a explicit formula for convolution-type asymmetric kernels:
$$
\lim_{l\to\infty}\lambda\!\left((-l,l)\right)
=\inf_{\nu\in\mathbb{R}}\left\{\int_{\mathbb{R}}J(x)e^{-\nu x}\,dx+c\nu\right\},
$$
expressing the whole-space principal eigenvalue as a one-dimensional variational problem involving the Laplace transform of the kernel. 
When the drift coefficient no longer has a fixed sign, even more striking phenomena emerge in the spectral structure: Ducrot, Griette, and Liang \cite{DucrotaGrietteLiang2026} proved in a periodically heterogeneous setting that, under sign-changing drift, the principal eigenproblem may admit a continuum of principal eigenpairs, all embedded in the continuous spectrum, with eigenvectors possessing singular measure components even appearing at the boundary --- a fundamental contrast to the isolated, simple principal eigenvalue of local diffusion operators. 
Most recently, Liang, Zhang, and Zhang \cite{LiangZhangZhang} extended the theory from one-dimensional/flat-torus settings to general compact Riemannian manifolds, establishing a necessary and sufficient criterion based on the spectral radius of the resolvent.

In this paper, we consider the generalized principal eigenvalues for the nonlocal operator
$$
\mathcal{L}_{\nu}^{\Omega}\phi = e^{-\nu x}\mathcal{L}^{\Omega}(e^{\nu \cdot}\phi)(x)
= \int_{\Omega} e^{-\nu(x-y)}J(x,y)\phi(y) \, dy + q(x)\phi'(x) + (a(x)+\nu q(x))\phi(x).
$$
From now on, when $\Omega = \mathbb R$, we shall omit the superscript and denote $\mathcal{L}^{\Omega}$ and $\mathcal{L}^{\Omega}_\nu$ simply by $\mathcal{L}$ and $\mathcal{L}_\nu$, respectively.
We mainly consider three types of admissible sets:
$$
\mathcal{A}_1 := \Big\{\phi\in C^1(\mathbb{R}): 0<\inf_{\mathbb{R}}\phi \le \sup_{\mathbb{R}}\phi < \infty,\ \|\phi'\|_{L^\infty}<\infty\Big\},
$$
$$
\mathcal{A}_2 := \Big\{\phi\in C(\mathbb{R}): 0<\inf_{\mathbb{R}}\phi \le \sup_{\mathbb{R}}\phi < \infty,\ \phi \text{ is uniformly continuous}\Big\},
$$
$$
\mathcal{A}_\Omega := \Big\{\phi\in C^1(\Omega)\cap C(\overline{\Omega}): \phi>0 \text{ on } \Omega\Big\}.
$$
In particular, \(\mathcal{A}_{\mathbb{R}} = \{\phi\in C^1(\mathbb{R}): \phi>0 \text{ on } \mathbb{R}\}\). When \(\mathcal{A}=\mathcal{A}_2\), the inequalities in (1.1) are understood in the viscosity sense.
	The generalized principal eigenvalues associated with operator \(\mathcal{L}^\Omega_\nu\) are:
	\begin{equation*}
		\left\{
		\begin{aligned}
		\underline{\lambda}(\nu,\mathcal{A}_\Omega) &:= \sup\{\lambda \mid \exists \phi\in\mathcal{A}_\Omega \text{ s.t. } \mathcal{L}^\Omega_\nu\phi \ge \lambda\phi \text{ in } \Omega\}, \\
		\overline{\lambda}(\nu,\mathcal{A}_\Omega) &:= \inf\{\lambda \mid \exists \phi\in\mathcal{A}_\Omega \text{ s.t. } \mathcal{L}^\Omega_\nu\phi \le \lambda\phi \text{ in } \Omega\}, 
		\end{aligned}
		\right. 
	\end{equation*} 
   \begin{equation*}
   	\left\{
   	\begin{aligned}
   		\underline{\lambda}(\nu,\mathcal{A}_1) &:= \sup\{\lambda \mid \exists \phi\in\mathcal{A}_1 \text{ s.t. } \mathcal{L}_\nu\phi \ge \lambda\phi \text{ in } \mathbb{R}\}, \\
   	\overline{\lambda}(\nu,\mathcal{A}_1) &:= \inf\{\lambda \mid \exists \phi\in\mathcal{A}_1 \text{ s.t. } \mathcal{L}_\nu\phi \le \lambda\phi \text{ in } \mathbb{R}\}, 
   		\end{aligned}
   	\right. 
   \end{equation*} 
   and
   \begin{equation*}
   	\left\{
   	\begin{aligned}
   			\underline{\lambda}(\nu,\mathcal{A}_2) &:= \sup\{\lambda \mid \exists \phi\in\mathcal{A}_2 \text{ s.t. } \mathcal{L}_\nu\phi \ge \lambda\phi \text{ in } \mathbb{R} \text{ in the viscosity sense}\}, \\
   		\overline{\lambda}(\nu,\mathcal{A}_2) &:= \inf\{\lambda \mid \exists \phi\in\mathcal{A}_2 \text{ s.t. } \mathcal{L}_\nu\phi \le \lambda\phi \text{ in } \mathbb{R} \text{ in the viscosity sense}\}.
   			\end{aligned}
   	\right. 
   \end{equation*}

\begin{definition}
	An upper semi-continuous (USC) function \(\phi:\mathbb{R}\to\mathbb{R}\) is said to satisfy
		$$
	\int_{\mathbb{R}} e^{-\nu(x-y)}J(x,y)\phi(y)\,dy + q(x)\phi'(x) + (a(x)+\nu q(x))\phi(x) \ge F(x)
	$$
 in the viscosity sense	if for any test function \(\Phi\in C^1\) such that
	\begin{enumerate}
		\item[(a)] \(\Phi(x_0)=\phi(x_0)\) for some \(x_0\in\mathbb{R}\),
		\item[(b)] \(\Phi(x)\ge \phi(x)\) for all \(x\in\{x_0+\mathcal{N}\}\setminus\{x_0\}\) with some neighborhood \(\mathcal{N}\) of the origin,
	\end{enumerate}
	then
	$$
	\int_{\mathbb{R}} e^{-\nu(x_0-y)}J(x_0,y)\phi(y)\,dy + q(x_0)\Phi'(x_0) + (a(x_0)+\nu b(x_0))\phi(x_0) \ge F(x_0).
	$$
	We also sometimes refer to $\phi$ as a viscosity subsolution of the equation 
	$$
	\int_{\mathbb{R}} e^{-\nu(x-y)}J(x,y)\phi(y)\,dy + q(x)\phi'(x) + (a(x)+\nu q(x))\phi(x) = F(x).$$
	A lower semi-continuous (LSC) function \(\phi:\mathbb{R}\to\mathbb{R}\) is said to be a viscosity supersolution if \(-\phi\) is a viscosity subsolution. A continuous function is said to be a viscosity solution if it is both a viscosity subsolution and a viscosity supersolution.
\end{definition}
We impose the following assumptions on the bounded kernel \(J\):\\
(J1) \(J(x,\cdot)\) is measurable for all \(x\in \mathbb{R}\); there exist constants \(\kappa_0,\delta_0>0\) and a continuous function \(J_0\) such that
$$
\kappa_0 \mathbf{1}_{(x-\delta_0,x+\delta_0)}(y) \leq J(x,y) \leq J_0(x-y), \quad \forall (x,y)\in \mathbb{R}\times\mathbb{R},
$$
where \(J_0\) satisfies \(\displaystyle \int_{\mathbb{R}} J_0(y) e^{-\nu y} \, dy < \infty\) for any \(\nu\in\mathbb{R}\), and \(J_0\) is monotone when \(|y|\gg 1\).\\
(J2)	 \(J(x,x-y)\) is uniformly continuous in \(x\), locally uniformly with respect to \(y\). That is, for any \(\varepsilon>0\) and $M>0$, there exists \(\delta>0\) such that for all \(x_1,x_2\) with \(|x_1-x_2|<\delta\),
$$\sup\limits_{|y|\le M}|J(x_1,x_1-y)-J(x_2,x_2-y)|<\varepsilon.$$

We also assume that \(a\) and \(q\) satisfy\\
(A1) \(a,q \in C(\mathbb{R}) \cap L^\infty(\mathbb{R})\) and \(q\neq 0\) in \(\mathbb{R}\).\\
(A2) \(\inf_{\Omega} |q| > 0\).\\
(A3) \(q\) is globally Lipschitz continuous on \(\mathbb{R}\) with Lipschitz constant $L_q$.

Occasionally, we impose stronger assumptions than (J2) and (A1) as follows:\\
(J2)$^\prime$  \(J(x,x-y)\) is almost periodic in \(x\), locally uniformly with respect to \(y\). That is, for any \(\varepsilon>0\), $M>0$, and sequence \(\{x_n\}\), there is a subsequence (still denoted by \(\{x_n\}\)) and $n_0$ such that
$$
\sup\limits_{x\in\mathbb R,\ |y|\le M}\big| J(x + x_m, x + x_m - y) - J(x + x_n, x + x_n - y) \big|  \le \varepsilon \quad \text{ for any } m, n \ge n_0.
$$
%uniformly with respect to \(x \in \mathbb{R}\).\\
(A1)$^\prime$ $a,q$ are almost periodic functions.

    The main results of this paper are as follows.
\begin{theorem}\label{thm:main_result1}
	The following are valid:\\
	(i) It always holds that $\underline{\lambda}(\nu,\mathcal{A}_1) \le \underline{\lambda}(\nu,\mathcal{A}_2)$ and $\overline{\lambda}(\nu,\mathcal{A}_2) \le \overline{\lambda}(\nu,\mathcal{A}_1) \text{ for any }\nu\in\mathbb R$.\\
	(ii) Assume (J1) holds and that \(q(x)\) and \(a(x)\) are bounded. Then
	$$
	\underline{\lambda}(\nu,\mathcal{A}_1) \le\overline{\lambda}(\nu,\mathcal{A}_1) \text{ for any }\nu\in\mathbb R.
	$$
	In addition, if (J2) and (A3) hold, and $a$ is uniformly continuous, then
		\begin{equation}\label{eq:GPE between A1 and A2}
				\underline{\lambda}(\nu,\mathcal{A}_1) \le \underline{\lambda}(\nu,\mathcal{A}_2) \le \overline{\lambda}(\nu,\mathcal{A}_2) \le \overline{\lambda}(\nu,\mathcal{A}_1) \text{ for any }\nu\in\mathbb R.
		\end{equation}
\end{theorem}

Next, we consider the generalized principal eigenvalues defined with respect to the test function set \(\mathcal{A}_{\Omega}\). 
First, it is easy to see from the definition that, for any \(\Omega\subset\mathbb{R}\),
\(\overline{\lambda}(\nu,\mathcal{A}_{\Omega})\) and \(\underline{\lambda}(\nu,\mathcal{A}_{\Omega})\) are independent of \(\nu\).

\begin{theorem}\label{thm:main_result2}
	Assume that \(J,q,a\) satisfy (J1), (J2), (A1), and (A2). Then the following statements are valid:\\
	(i) If \(\Omega=(r_1,r_2)\subset\mathbb{R}\) is bounded, then
		$$
		\overline{\lambda}(0,\mathcal{A}_\Omega)=
		\begin{cases}
			\sup\{\lambda \mid \exists \phi\in\mathcal{A}_\Omega,\ \phi(r_1)=0,\ \mathcal{L}^\Omega\phi\ge \lambda\phi \text{ on } \Omega\}, & q<0, \\[6pt]
			\sup\{\lambda \mid \exists \phi\in\mathcal{A}_\Omega,\ \phi(r_2)=0,\ \mathcal{L}^\Omega\phi\ge \lambda\phi \text{ on } \Omega\}, & q>0,
		\end{cases}
		$$
		and there exists \(\phi\in\mathcal{A}_\Omega\) such that \(\mathcal{L}\phi=\overline{\lambda}(0,\mathcal{A}_\Omega)\phi\), with \(\phi(r_1)=0\) if \(q<0\), \(\phi(r_2)=0\) if \(q>0\).\\		
	(ii) For any sequence \(\{\Omega_n\}_{n\in\mathbb{Z}_+}\) of non-empty bounded domains such that \(\overline{\Omega}_n\subset\Omega_{n+1}\subset\Omega\) for all \(n\in\mathbb{N}\) and \(\bigcup_{n\in\mathbb{N}}\Omega_n=\Omega\), one has
		$$
		\overline{\lambda}(0,\mathcal{A}_{\Omega_n}) \nearrow   \overline{\lambda}(0,\mathcal{A}_\Omega), \text{ as }n\to\infty
		$$
		and there exists \(\phi\in\mathcal{A}_{\Omega}\) such that \(\mathcal{L}^\Omega\phi=\overline{\lambda}(0,\mathcal{A}_{\Omega})\phi$.\\
    (iii) For any domain $\Omega$, whether bounded or unbounded, one has $\underline{\lambda}(0,\mathcal{A}_\Omega)=+\infty$. Therefore,
		$		\overline{\lambda}(0,\mathcal{A}_\Omega) < \underline{\lambda}(0,\mathcal{A}_\Omega).$
\end{theorem}
When the kernel has compact support, the above result can be found in \cite{Coville2020}.
The following theorem examines the relationship between
\(\underline{\lambda}(\nu,\mathcal{A}_{1})\),
\(\overline{\lambda}(\nu,\mathcal{A}_{1})\) and
\(\overline{\lambda}(0,\mathcal{A}_{\mathbb{R}})\).
To consider a slightly more general situation, we assume here that\\
(J1)$^\prime$ \(J(x,\cdot)\) is measurable for all \(x\in \mathbb{R}\); there exist constants \(\kappa_0,\delta_0>0\) such that
$$
\kappa_0 \mathbf{1}_{(x-\delta_0,x+\delta_0)}(y) \leq J(x,y), \quad \forall (x,y)\in \mathbb{R}\times\mathbb{R},
$$
and assume that $$\mathcal{I}_x=\left\{\nu\in\mathbb{R}:\int_{\mathbb{R}}J(x,x-y)e^{-\nu y}dy<+\infty\right\},$$
 is an interval independent of \(x\), which is denoted by \(\mathcal{I}\). 
 
When \(\nu\notin\mathcal{I}\), we adopt the convention that
$$
\underline{\lambda}(\nu,\mathcal{A}_{1})
=\overline{\lambda}(\nu,\mathcal{A}_{1})
=+\infty.
$$
%Then \((-\eta_0,\eta_0)\subset\mathcal{I}\equiv \mathcal I_x\) for some \(\eta_0>0\).
\begin{theorem}\label{the:combining hamel and du}
	Assume that  \(J,q,a\) satisfy (J1)$^\prime$, (J2), (A1), and (A2).  Then 
   \begin{equation}\label{eq:1 of the:combining hamel and du}
   	\inf\limits_{\nu\in\mathbb R}\underline{\lambda}(\nu,\mathcal A_1)\leq\lim\limits_{n\to\infty}\overline{\lambda}(0,\mathcal{A}_{\Omega_n})\leq\overline\lambda(0,\mathcal A_\mathbb R)\le\inf\limits_{\nu\in\mathbb R}\overline\lambda(\nu,\mathcal A_1),
   \end{equation}
	where $\{\Omega_n\}_n$ was given by (ii) of Theorem \ref{thm:main_result2} with $\Omega=\mathbb R.$
In particular, if \(\mathcal{I} = \{0\}\), i.e., \(J\) is heavy-tailed, then (3) becomes
$$
\underline{\lambda}(0, \mathcal{A}_1) \le \lim_{n\to\infty} \overline{\lambda}(0, \mathcal{A}_{\Omega_n}) \le \overline{\lambda}(0, \mathcal{A}_\mathbb{R}) \le \overline{\lambda}(0, \mathcal{A}_1).
$$
\end{theorem}
\begin{remark} We emphasize the connections of Theorem \ref{the:combining hamel and du} with the previous literature:
	
	(i) {Connection with \cite{Coville2020}:} In \cite{Coville2020}, the domain approximation property \(\lim_{n\to\infty}\overline{\lambda}(0,\mathcal{A}_{\Omega_n}) = \overline{\lambda}(0,\mathcal{A}_{\mathbb{R}})\) was established for kernels with compact support.
	 We generalize this result in Theorem 1.3(ii) to kernels with super-exponential decay (see (J1)).
	  Furthermore, in the case where the dispersal kernel requires no global upper bound (see (J1)'), inequality (3) also provides a comprehensive two-sided bound for the domain approximation limit via the exponential variational parameter \(\nu\).
	
	(ii) {Connection with \cite{DuFangNi2025}:} For convolution-type asymmetric kernels,  \cite{DuFangNi2025} obtained the explicit variational formula for the whole-space principal eigenvalue: \(\lim_{l\to\infty}\lambda((-l,l)) = \inf_{\nu \in \mathbb{R}}\left\{\int_{\mathbb{R}}J(y)e^{-\nu y}dy + \nu c\right\}\). 
	Our Theorem \ref{the:combining hamel and du} extends this variational characterization from convolution kernels to spatially heterogeneous kernels \(J(x,y)\) and spatially varying coefficients \(a(x), q(x)\).
	
	(iii) {Equality case:} Under the conditions of Theorem \ref{the:combining hamel and du}, if there exists \(\Phi \in \mathcal{A}_1\) such that \(\mathcal{L}_{\nu_0}\Phi = \underline\lambda(\nu_0, \mathcal{A}_1)\Phi\), where \(\nu_0\) satisfies \(\underline\lambda(\nu_0, \mathcal{A}_1) = \inf_{\nu \in \mathbb{R}} \underline\lambda(\nu, \mathcal{A}_1)\), then it is easy to see that \(\underline\lambda(\nu_0, \mathcal{A}_1) \ge \overline{\lambda}(0, \mathcal{A}_{\mathbb{R}})\). 
	Consequently, from the above theorem, we have
	$$
	\lim_{n\to\infty}\overline{\lambda}(0, \mathcal{A}_{\Omega_n}) = \overline{\lambda}(0, \mathcal{A}_{\mathbb{R}}) = \inf_{\nu \in \mathbb{R}} \underline\lambda(\nu, \mathcal{A}_1).
	$$
	Moreover, if \(\underline{\lambda}(\nu,\mathcal{A}_1) = \overline{\lambda}(\nu,\mathcal{A}_1)\) for any \(\nu \in \mathbb{R}\) (which holds, for instance, if \(\mathcal{L}_\nu\) possesses a genuine bounded principal eigenfunction bounded away from zero), the inequalities in (3) collapse into an equality.
	As a direct verification, consider the special case in  \cite{DuFangNi2025} where \(J(x,y) = J(x-y)\), \(a = 0\), and \(q = c\) is a constant. Then \(\mathcal{L}_\nu\) has the eigenpair \(\left(\int_{\mathbb{R}}e^{-\nu y}J(y)dy + \nu c, 1\right)\), and the equality above exactly recovers the explicit formula in  \cite{DuFangNi2025}. 
	Thus, Theorem \ref{the:combining hamel and du} bridges the domain approximation result in \cite{Coville2020} and the variational formula in  \cite{DuFangNi2025} under a unified, less restrictive framework.
\end{remark}
In what follows, we prove that \(\underline{\lambda}(\nu,\mathcal{A}_i)=\overline{\lambda}(\nu,\mathcal{A}_i)\) holds under almost periodic coefficients for \(\nu\in\mathbb{R}\) and \(i=1,2\).
\begin{theorem}\label{thm:main_result3}
Assume that \(J, q, a\) satisfy (J1), (J2)\(^\prime\), (A1)\(^\prime\), (A2), and (A3). 
Then
$$
\underline{\lambda}(\nu,\mathcal{A}_1) = \underline{\lambda}(\nu,\mathcal{A}_2) = \overline{\lambda}(\nu,\mathcal{A}_2) = \overline{\lambda}(\nu,\mathcal{A}_1) \ \forall \nu\in\mathbb R.
$$
We denote this generalized eigenvalue simply by $\lambda(\nu)$. Then for any sequence \(\{\Omega_n\}_{n\in\mathbb{Z}_+}\) of non-empty bounded domains such that \(\overline{\Omega}_n\subset\Omega_{n+1}\) for all \(n\in\mathbb{N}\) and \(\bigcup_{n\in\mathbb{N}}\Omega_n=\mathbb{R}\), we have
$$
\overline{\lambda}(0,\mathcal{A}_{\mathbb{R}}) =\lim_{n\to\infty}\overline{\lambda}(0,\mathcal{A}_{\Omega_n}) = \inf_{\nu\in\mathbb{R}}{\lambda}(\nu).
$$
\end{theorem}

\textbf{Outline of this paper:}
The rest of this paper is organized as follows. In Section 2, we introduce the main tool of this paper: a Harnack type inequality for the nonlocal operator. Section 3 is devoted to the general coefficients case, where we provide detailed proofs of Theorems \ref{thm:main_result1}, \ref{thm:main_result2} and \ref{the:combining hamel and du}. In Section 4, we focus on the almost periodic coefficients case and give the proof of Theorem \ref{thm:main_result3}.

\begin{comment}
\begin{definition}
	Consider the operator \(\mathcal{L}\) on \(\mathbb{R}\). We say \(\mathcal{L}\) is
	\begin{enumerate}
		\item[(i)] supercritical if the set \(\{\phi \mid \phi>0 \text{ on } \mathbb{R},\ \mathcal{L}\phi=0\}\) is empty;
		\item[(ii)] critical if the set \(\{\phi \mid \phi>0 \text{ on } \mathbb{R},\ \mathcal{L}\phi\le 0\}\) is one-dimensional;
		\item[(iii)] subcritical if it is neither critical nor supercritical.
	\end{enumerate}
\end{definition}
\begin{theorem}\label{thm:main_result4}
	Assume (J1) holds, \(q(x)\) and \(a(x)\) are bounded. Then
	$$
	(\mathcal{L}-\lambda) \text{ is supercritical} \iff \lambda < \overline{\lambda}(0,\mathcal{A}_{\mathbb{R}}),
	$$
	$$
	(\mathcal{L}-\lambda) \text{ is critical or supercritical if } \lambda = \overline{\lambda}(0,\mathcal{A}_{\mathbb{R}}),
	$$
	$$
	(\mathcal{L}-\lambda) \text{ is subcritical if } \lambda > \overline{\lambda}(0,\mathcal{A}_{\mathbb{R}}).
	$$
\end{theorem}
\end{comment}

\section{The main tool: Harnack type inequality}
  In this section, we generalize the Harnack-type inequality established in \cite[Lemma 3.1]{Coville2020} to kernels that are not necessarily compactly supported.
   This generalized estimate will play a key role in our investigation of the relations between the various notions of generalized principal eigenvalues considered in this paper.
\begin{theorem}\label{harnack}
	Let \(\Omega\subseteq\mathbb{R}\) be a non-empty domain. Assume that \(J\) satisfies (J1), (J2), and \(q\) satisfies (A1), (A2). Then for any \(r_1<r_2\) such that \([r_1,r_2]\subset\Omega\), there exists a positive constant \(C\) depending only on \(|r_2-r_1|\), \(\operatorname{dist}([r_1,r_2],\partial\Omega)\), \(J\), \(J_0\), \(q\), and \(a\) such that
	$$
	\max_{[r_1,r_2]} u \le C \min_{[r_1,r_2]} u,
	$$
	where \(u\in C^1(\Omega)\cap C(\overline{\Omega})\) is a positive solution to
	\begin{equation}\label{eq:harnack1}
		\int_{\Omega} J(x,y)u(y)\,dy + q(x)u'(x) + a(x)u(x) = 0,\quad x\in\Omega. 
	\end{equation}
\end{theorem}

\begin{proof}
	The proof is divided into five steps. 
	Denote  $\Omega$ by $(\alpha,\beta)$.
	
	\medskip
	\noindent{Step 1: A preliminary simplification.}
	
	Without loss of generality, we may assume that \(q>0\). In fact, if \(q<0\), set \(\widetilde{u}(x)=u(-x)\), \(\widetilde{q}(x)=-q(-x)\), \(\widetilde{a}(x)=a(-x)\), \(\widetilde{J}(x,y)=J(-x,-y)\), \(x,y\in\widetilde{\Omega}:=\{x: -x\in\Omega\}\). One has
	$$
	\int_{\widetilde{\Omega}}\widetilde{J}(x,y)\widetilde{u}(y)\,dy + \widetilde{q}(x)\widetilde{u}'(x) + \widetilde{a}(x)\widetilde{u}(x)=0,\quad x\in\widetilde{\Omega}.
	$$
	Note that \(\inf_{x\in\widetilde{\Omega}}\widetilde{q}(x)=q_0>0\). Then one can still deduce the Harnack type inequality.
	
	Next, motivated by \cite{Coville2020} Lemma 3.1, define
	$$
	w(x):=e^{\int_{r_1}^x (Q(s)+1)ds}u(x),\qquad \widetilde{w}(x):=e^{\int_{r_1}^x Q(s)ds}u(x)=e^{-x+r_1}w(x),
	$$
	$$
	K(x,y):=\frac{J(x,y)}{q(x)}e^{\int_y^x (Q(s)+1)ds},\qquad \widetilde{K}(x,y):=\frac{J(x,y)}{q(x)}e^{\int_y^x Q(s)ds},
	$$
	where \(Q(s):=a(s)/q(s)\). From \eqref{eq:harnack1} and (A1), \(w\) and \(\widetilde{w}\) are positive in \(\Omega\), of class \(C^1(\Omega)\cap C(\overline{\Omega})\), and satisfy
	\begin{equation}\label{eq:harnack2}
		w'(x)+\int_{\Omega}K(x,y)w(y)\,dy - w(x)=0,\quad x\in\Omega, 
	\end{equation}
	\begin{equation}\label{eq:harnack3}
		\widetilde{w}'(x)+\int_{\Omega}\widetilde{K}(x,y)\widetilde{w}(y)\,dy=0,\quad x\in\Omega. 
	\end{equation}
	
	\(K\) and \(\widetilde{K}\) satisfy (J1) and (J2). If we can prove a Harnack type inequality for \(w\) solving \eqref{eq:harnack2}, i.e.,
	$$
	\max_{[r_1,r_2]} w(x) \le C \min_{[r_1,r_2]} w(x),
	$$
	then we can obtain a Harnack type estimate for \(u\) using \(w(x)=e^{\int_{r_1}^x(Q(s)+1)ds}u(x)\).
	
	\medskip
	\noindent{Step 2: Proof of a special case.}
	
	Let \(x_0\) be a point where \(w\) achieves its maximum in \([r_1,r_2]\). Then either \(x_0=r_2\), or \(x_0\in[r_1,r_2)\) and \(w'(x_0)\le 0\). Note that \(\widetilde{w}\) is decreasing by \eqref{eq:harnack3}. If \(x_0=r_2\), then monotonicity of \(\widetilde{w}\) yields
	\begin{equation*}
		\begin{split}
		&\quad	w(x_1)\le w(x_0)=w(r_2)=\widetilde{w}(r_2)e^{r_2-r_1}\\
		&\le \widetilde{w}(x_2)e^{r_2-r_1}=w(x_2)e^{-x_2+r_1}e^{r_2-r_1}\le e^{r_2-r_1}w(x_2),
		\end{split}
	\end{equation*}
	for any \(x_1,x_2\in[r_1,r_2]\). This gives the desired Harnack type inequality.
	
	Hence we only need to prove a Harnack type inequality for \(w\) under the conditions \(q_0=\inf_{\Omega}q>0\) and \(x_0\in[r_1,r_2]\) with \(w'(x_0)\le 0\).
	
	To do so, we first claim that if \(\alpha>-\infty\), then there exists \(d^->0\) small such that
	\begin{equation}
		\widetilde{w}(x)\ge \frac12 \widetilde{w}(\alpha)\quad \text{for } x\in[\alpha,\alpha+2d^-]. 
	\end{equation}
	
	In fact, it follows from \eqref{eq:harnack3} that
	$$
	\widetilde{w}'(x)=-\int_{\alpha}^{\beta}K(x,y)\widetilde{w}(y)\,dy
	\ge -\int_{\alpha}^x \widetilde{K}(x,y)\,dy\,\widetilde{w}(\alpha)-\int_x^\beta \widetilde{K}(x,y)\,dy\,\widetilde{w}(x).
	$$
	Let $\|Q\|_\infty:=\sup\limits_{s\in\mathbb R}|a(s)/q(s)|$. Denote \(\frac{\int_0^{+\infty}J(y)e^{\|Q\|_\infty y}dy}{q_0}\) and \(\frac{\int_{-\infty}^0 J(y)e^{\|Q\|_\infty y}dy}{q_0}\) by \(J_0^+\) and \(J_0^-\), respectively. Then
	$$
	\widetilde{w}(x)\ge \widetilde{w}(\alpha)\left[e^{-J_0^+(x-\alpha)}-\frac{J_0^-}{J_0^+}\left(1-e^{-J_0^+(x-\alpha)}\right)\right].
	$$
	One can therefore obtain the desired \(d^-\) from the above equation.
	
	Set
	$$
	d=\min\left\{d^-,\frac{r_2-r_1}{2},\frac{r_1-\alpha}{2},\frac{\beta-r_2}{2},1\right\}.
	$$
	
	\medskip
	\noindent{Step 3: Proof of the Harnack type inequality  when \(-\infty<\alpha<\beta<+\infty\).}
	
	It follows from \eqref{eq:harnack2} that
	\begin{equation}\label{eq:estimate on w(x0)}
		\begin{split}
		&\quad	w(x_0) = w'(x_0)+\int_{\alpha}^{\beta}K(x_0,y)w(y)\,dy\le \int_{\alpha}^{\beta}K(x_0,y)w(y)\,dy  \\
			&= \underbrace{\int_{\alpha}^{\alpha+d}K(x_0,y)w(y)\,dy}_{I_1}
			+\underbrace{\int_{\beta-d}^{\beta}K(x_0,y)w(y)\,dy}_{I_2}
			+\underbrace{\int_{\alpha+d}^{\beta-d}K(x_0,y)w(y)\,dy}_{I_3}. 
		\end{split}
	\end{equation}
	
	Note that
	\begin{equation}\label{eq:estimate on w(x0)1}
		\begin{split}
			I_1 &= \int_{\alpha}^{\alpha+d}K(x_0,y)e^{y-r_1}\widetilde{w}(y)\,dy \\
			&\le \int_{\alpha}^{\alpha+d}K(x_0,y)e^{y-r_1}\widetilde{w}(\alpha)\,dy \\
			&\le 2e^{-d}\int_{\alpha+d}^{\alpha+2d}K(x_0,y-d)w(y)\,dy  \\
			&\le \frac{2e^{-d}}{q_0}\sup\limits_{\xi\in[r_1-\alpha-d,r_2-\alpha]}J_0(\xi)e^{(\|Q\|_\infty+1)\xi}\int_{\alpha+d}^{\alpha+2d}w(y)\,dy, 
		\end{split}
	\end{equation}
	and
	\begin{equation}\label{eq:estimate on w(x0)2}
		I_2 \le \frac{e^d}{q_0}\sup_{\xi\in[r_1-\beta+d,r_2-\beta]}J_0(\xi)e^{(\|Q\|_\infty+1)\xi}\int_{\beta-2d}^{\beta-d}w(y)\,dy,
	\end{equation}
	and
	\begin{equation}\label{eq:estimate on w(x0)3}
		I_3 \le \frac{1}{q_0}\sup_{\xi\in[r_1+d-\beta,r_2-d-\alpha]}J_0(\xi)e^{(\|Q\|_\infty+1)\xi}\int_{\alpha+d}^{\beta-d}w(y)\,dy. 
	\end{equation}
	
	Combining \eqref{eq:estimate on w(x0)}-\eqref{eq:estimate on w(x0)3}, we have
	\begin{equation}\label{eq:w(x0) smaller than}
		w(x_0) \le \frac{1}{q_0}(2e^{-d}+e^d+1)\sup_{\xi\in[r_1+d-\beta,r_2-d-\alpha]}J_0(\xi)e^{(\|Q\|_\infty+1)\xi}\int_{\alpha+d}^{\beta-d}w(y)\,dy. 
	\end{equation}
	
	Denote \(\delta:=\min\{d_0/2,\delta/4\}\). Then it follows from Proposition \ref{prop:claim} below that there exists a positive constant $C_1>0$, depending only on $\operatorname{dist}([r_1,r_2],\partial\Omega)$, $J$, $q$, and $a$, such that
	\begin{equation}\label{eq:cliam}
	\int_{x_1-\delta/4}^{x_1+\delta/4}w(y)\,dy \ge C_1\int_{x_2-\delta/4}^{x_2+\delta/4}w(y)\,dy
	\quad \text{for all } x_1,x_2\in[\alpha+d,\beta-d].
	\end{equation}
	Taking any \(x_1\in[r_1,r_2]\), we have
	\begin{equation*}
		\begin{split}
		\widetilde{w}(x_1)&=\widetilde{w}(x_1+\delta)+\int_{x_1}^{x_1+\delta}\int_{\alpha}^{\beta}\widetilde{K}(x,y)\widetilde{w}(y)\,dy\,dx\\
		&\ge \delta\frac{\kappa_0}{\|q\|_\infty}e^{-\|Q\|_\infty\delta/4}\int_{x_1-\delta/4}^{x_1+\delta/4}\widetilde{w}(y)\,dy,
	\end{split}
	\end{equation*}
	which yields
	\begin{equation}\label{eq:w(x1)larger than}
		w(x_1)\ge \frac{\delta e^{-\delta}\kappa_0e^{-\|Q\|_\infty\delta/4}}{\|q\|_\infty}\int_{x_1-\delta/4}^{x_1+\delta/4}w(y)\,dy. 
	\end{equation}
	
	Furthermore, setting \(N:=\min\{n\in\mathbb{N}: n\ge 4(\beta-\alpha-2d)/\delta\}\), we have
	$$
	[\alpha+d,\beta-d]=\bigcup_{k=0}^{N-1}\left[\alpha+d+\frac{k\delta}{4},\alpha+d+\frac{(k+1)\delta}{4}\right].
	$$
	Then by \eqref{eq:cliam},
	\begin{equation}\label{eq:4.11}
		\int_{x_1-\delta/4}^{x_1+\delta/4}w(y)\,dy \ge \frac{C_1}{N}\int_{\alpha+d}^{\beta-d}w(y)\,dy. 
	\end{equation}
	
	Combining \eqref{eq:w(x0) smaller than}-\eqref{eq:4.11}, we obtain
	$$
	w(x_0)\le \widetilde{C}\cdot \frac{N}{C_1}\cdot \frac{\|q\|_\infty}{\delta\kappa_0}e^{\delta}e^{\|Q\|_\infty\delta/4}w(x_1).
	$$
	Thus completes the proof in this step.
	
	\medskip
	\noindent{Step 4: Proof of the Harnack type inequality  when \(-\infty<\alpha<\beta=+\infty\).}
	
	It follows from \eqref{eq:harnack2} that
    \begin{equation*}
    	\begin{split}
    		&\quad	w(x_0)\le \int_{\alpha}^{+\infty}K(x_0,y)w(y)\,dy\\
    		&= I_1+\int_{\alpha+d}^{x_0+\frac{N_0+1}{2}\delta}K(x_0,y)w(y)\,dy+\int_{x_0+\frac{N_0+1}{2}\delta}^{+\infty}K(x_0,y)w(y)\,dy,
    	\end{split}
	\end{equation*}
	where \(N_0\) is large such that \(J_0\) is increasing on \((-\infty,-N_0\delta/2]\).
	
	Similar estimates as in Step 3 give
	$$
	w(x_0)\le C\int_{\alpha+d}^{r_2+N_0\delta}w(y)\,dy.
	$$
	Then, similarly to Step 3, one obtains
	$$
	w(x_0)\le \widetilde{C}w(x_1)\quad \text{for all } x_1\in[r_1,r_2].
	$$
	
	\medskip
	\noindent{Step 5: Proof of the Harnack type inequality when \(-\infty=\alpha<\beta<+\infty\) or \((\alpha,\beta)=\mathbb{R}\).}
	
	If \(-\infty=\alpha<\beta<+\infty\), then from \eqref{eq:harnack2},
	$$
	w(x_0)\le \int_{-\infty}^{\beta}K(x_0,y)w(y)\,dy
	= \sum_{n=N_1+1}^{\infty}\int_{\widetilde{T}_n}K(x_0,y)w(y)\,dy
	+\int_{x_0-\frac{N_1+1}{2}\delta}^{\beta-d}K(x_0,y)w(y)\,dy+I_3,
	$$
	where \(N_1\) is large such that \(J_0\) is decreasing on \([N_0\delta/2,+\infty)\), and \(\widetilde{T}_n:=[x_0-(n+1)\delta/2,x_0-n\delta/2]\). Similar estimates yield
	$$
	w(x_0)\le \widetilde{C}_3\int_{r_1-\frac{N_1+1}{2}\delta}^{\beta-d}w(y)\,dy,
	$$
	which is the analogue of \eqref{eq:w(x0) smaller than}. Thus a Harnack type inequality follows.
	
	If \((\alpha,\beta)=\mathbb{R}\), then
	$$
	w(x_0)\le \left(\widetilde{C}_1+\widetilde{C}_4+\frac{2\max\{1,e^{-\nu\delta}\}}{\delta q_0}\int_0^{+\infty}J_0(\xi)e^{\nu\xi}\,d\xi\right)\int_{r_1-\frac{N_1+1}{2}\delta}^{r_2+\frac{N_0+1}{2}\delta}w(y)\,dy,
	$$
	which is the analogue of \eqref{eq:w(x0) smaller than} in Step 3. One thus obtains a Harnack type inequality in this case as before.
\end{proof}
\begin{remark}\label{re:more-about-HIE}
    (i) Consider any family of domains \( \Omega_s \subset \mathbb{R} \) indexed by \( s \in \Lambda \). Let \( u_s \) be a solution of
     $$
     \int_{\Omega_s} J(x,y) u_s(y) \, dy + q_s(x) u_s'(x) + a(x) u_s(x) = 0, \quad x \in \Omega_s.
     $$
     If \( [r_1, r_2] \subset \Omega_s \) and
     $$
     \inf_{s \in \Lambda} \operatorname{dist}([r_1, r_2], \partial \Omega_s) \ge d > 0,
     $$
     then by the proof of  Theorem \ref{harnack}, there exists a constant \( C > 0 \), depending only on \( |r_2 - r_1| \), \( d \), \( J \), \( J_0 \), \( q \), and \( a \), such that for all \( s \in \Lambda \),
     $$
     \max_{[r_1, r_2]} u_s \le C \min_{[r_1, r_2]} u_s.
     $$	
     (ii) From the proof, it is evident that if \(\nu_0\) is sufficiently large (e.g., \(\nu_0 > \|Q\|_\infty + 1\)), and \(J_0\) satisfies \(\displaystyle \int_{\mathbb{R}} J_0(y) e^{-\nu y} \, dy < \infty\) for some \(\nu \ge \nu_0\), then the Harnack inequality still holds.
\end{remark}
    We now establish the proposition that will be used in the proof of the above theorem. This proposition is essentially \cite[(3.13)]{Coville2020}; we include its proof here for the sake of completeness.    
\begin{proposition}\label{prop:claim}
	Let \(w(x)\) be the positive solution to \eqref{eq:harnack2} on \(\Omega\) (non-empty bounded or unbounded domain), where
	$$
	\delta\in\left(0,\min\left\{\frac{\delta_0}{2},\frac{1}{2}\operatorname{dist}\big([\alpha_1,\beta_1],\partial\Omega\big)\right\}\right).
	$$
	Then there exists \(C_1>0\) depending only on \(\operatorname{dist}([\alpha_1,\beta_1],\partial\Omega)\), \(\delta\), \(J\), \(q\), and \(a\), such that
	\begin{equation}
		\int_{x_1-\delta/4}^{x_1+\delta/4}w(y)\,dy \ge C_1\int_{x_2-\delta/4}^{x_2+\delta/4}w(y)\,dy
		\quad \forall\,x_1,x_2\in[\alpha_1,\beta_1]\subset\Omega. \label{eq:A1}
	\end{equation}
\end{proposition}

\begin{proof}
	For any \(z\in[\alpha_1-\delta,\beta_1+\delta]\), integrating \eqref{eq:harnack3} from \(z\) to \(z+\delta\), we have
	\begin{equation}\label{eq:A2}
		\begin{split}
			\widetilde{w}(z) &= \widetilde{w}(z+\delta) + \int_{z}^{z+\delta}\int_{\Omega} \widetilde{K}(x,y)\widetilde{w}(y)\,dy\,dx \\
			&\ge \int_{\Omega}\widetilde{w}(y)\int_{z}^{z+\delta}\widetilde{K}(x,y)\,dx\,dy  \\
			&\ge \int_{z-\delta}^{z+\delta}\widetilde{w}(y)\int_{z}^{z+\delta}\widetilde{K}(x,y)\,dx\,dy \\
			&\ge \frac{\kappa_0\delta}{\|q\|_\infty}e^{-2\|Q\|_\infty\delta}\int_{z-\delta}^{z+\delta}\widetilde{w}(y)\,dy. 
		\end{split}
	\end{equation}
	Set
	$$
	\widetilde{\kappa}_0:=\frac{\kappa_0}{\|q\|_\infty}e^{-2\|Q\|_\infty\delta}.
	$$
	From \(\widetilde{w}(x)=e^{-x+r_1}w(x)\), we have
	\begin{equation}
		\widetilde{w}(z) \ge \widetilde{\kappa}_0\delta e^{-\delta}\int_{z-\delta}^{z+\delta}\widetilde{w}(y)\,dy. \label{eq:A3}
	\end{equation}
	
	Integrating \eqref{eq:A3} over \(\left[x_1-\frac{\delta}{4},x_1+\frac{\delta}{4}\right]\) with \(x_1\in[\alpha_1,\beta_1]\), we have
	\begin{align}\label{eq:A4}
		\begin{split}
			&\quad\int_{x_1-\delta/4}^{x_1+\delta/4}w(s)\,ds\ge \widetilde{\kappa}_0\delta e^{-\delta}\int_{x_1-\delta/4}^{x_1+\delta/4}\int_{s-\delta}^{s+\delta}w(y)\,dy\,ds  \\
			&\ge \widetilde{\kappa}_0\delta e^{-\delta}\left(
			\int_{x_1-\delta/4}^{x_1-\delta/8}\int_{s-\delta}^{s+\delta}w(y)\,dy\,ds
			+\int_{x_1+\delta/8}^{x_1+\delta/4}\int_{s-\delta}^{s+\delta}w(y)\,dy\,ds
			\right)  \\
			&\ge \widetilde{\kappa}_0\delta e^{-\delta}\left(
			\int_{x_1-\delta/4}^{x_1-\delta/8}\int_{x_1-\delta/2}^{x_1}w(y)\,dy\,ds
			+\int_{x_1+\delta/8}^{x_1+\delta/4}\int_{x_1}^{x_1+\delta/2}w(y)\,dy\,ds
			\right)  \\
			&= \frac{\widetilde{\kappa}_0\delta^2 e^{-\delta}}{8}\int_{x_1-\delta/2}^{x_1+\delta/2}w(y)\,dy. 
		\end{split}
	\end{align}
	Hence, for any \(x_1\in[\alpha_1,\beta_1]\),
	\begin{equation}\label{eq:A5}
		\int_{x_1-\delta/4}^{x_1+\delta/4}\widetilde{w}(s)\,ds
		\ge \frac{\widetilde{\kappa}_0\delta^2 e^{-\delta}}{8}\int_{z-\delta/4}^{z+\delta/4}\widetilde{w}(y)\,dy,
		\quad \forall\,z\in\left[x_1-\frac{\delta}{4},x_1+\frac{\delta}{4}\right]. 
	\end{equation}
	Note that \eqref{eq:A4} holds for any \(z\in\left[x_1-\frac{\delta}{4},x_1+\frac{\delta}{4}\right]\cap[\alpha_1,\beta_1]\). Then, \eqref{eq:A5} implies that
	\begin{equation}
		\int_{x_1-\delta/4}^{x_1+\delta/4}\widetilde{w}(s)\,ds
		\ge \left(\frac{\widetilde{\kappa}_0\delta^2 e^{-\delta}}{8}\right)^2
		\int_{z-\delta/2}^{z+\delta/2}\widetilde{w}(y)\,dy,
		\quad \forall\,z\in\left[x_1-\frac{\delta}{4},x_1+\frac{\delta}{4}\right]\cap[\alpha_1,\beta_1], \label{eq:A6}
	\end{equation}
	which yields that
	\begin{equation}
		\int_{x_1-\delta/4}^{x_1+\delta/4}\widetilde{w}(s)\,ds
		\ge \left(\frac{\widetilde{\kappa}_0\delta^2 e^{-\delta}}{8}\right)^2
		\int_{z-\delta/4}^{z+\delta/4}\widetilde{w}(y)\,dy,
		\quad \forall\,z\in\left[x_1-\frac{\delta}{2},x_1+\frac{\delta}{2}\right]\cap[\alpha_1,\beta_1]. \label{eq:A7}
	\end{equation}
	
	Taking $	N=\min\left\{n\in\mathbb{N}\,\bigg|\,n\ge \frac{4(\beta_1-\alpha_1)}{\delta}\right\}$
	and arguing by induction, we see that, for every \(k\in\{1,2,\dots,N\}\), and \(z\in\left[x_1-\frac{k\delta}{4},x_1+\frac{k\delta}{4}\right]\cap[\alpha_1,\beta_1]\),
	\begin{equation}
		\int_{x_1-\delta/4}^{x_1+\delta/4}w(s)\,ds\ge \left(\frac{\widetilde{\kappa}_0\delta^2 e^{-\delta}}{8}\right)^k \int_{z-\delta/4}^{z+\delta/4}w(y)\,dy. \label{eq:A8}
	\end{equation}	
	In particular,
	\begin{equation}
		\int_{x_1-\delta/4}^{x_1+\delta/4}w(s)\,ds\ge \left(\frac{\widetilde{\kappa}_0\delta^2 e^{-\delta}}{8}\right)^N \int_{z-\delta/4}^{z+\delta/4}w(y)\,dy \label{eq:A9}
	\end{equation}
	for any
	$$
	z\in\left[x_1-\frac{N\delta}{4},x_1+\frac{N\delta}{4}\right]\cap[\alpha_1,\beta_1]=[\alpha_1,\beta_1].
	$$
	Thus, we complete the proof by taking	$	C_1=\left(\frac{\widetilde{\kappa}_0\delta^2 e^{-\delta}}{8}\right)^N>0.$
\end{proof}

 \section{General coefficients: proofs of Theorems \ref{thm:main_result1} and \ref{thm:main_result2}}

 \subsection{Proof of Theorem \ref{thm:main_result1}}
 
Since the operator involved in Theorem 1.2 is defined on \(\Omega=\mathbb R\), we shall follow the convention introduced earlier and write \(\mathcal L\) for \(\mathcal L^\Omega\).
 To prove Theorem \ref{thm:main_result1}, we need the following lemma.
 
 \begin{lemma}\label{lem2.1}
   Assume that (J1) holds and that \(q(x)\) and \(a(x)\) are bounded. If \(z\ge 0\) satisfies, in the viscosity sense,
   $$
   \mathcal{L}z(x) \le -\varepsilon_0,\quad x\in\mathbb{R},
   $$
   for some \(\varepsilon_0>0\).
   Then \(\inf_{\mathbb{R}} z > 0\).
 \end{lemma}
 
 \begin{proof}
 	Suppose by contradiction that \(\inf_{x\in\mathbb{R}} z = 0\). Consider a nonnegative function \(\theta\in C^1(\mathbb{R})\) with \(\theta(0)=0\) and \(\lim_{|x|\to\infty}\theta(x)=1\). One can find a constant \(C\) large enough such that
 	$$
 	\mathcal{L}(\theta(x-\xi)) \le \frac{C}{4}\varepsilon_0,\quad \forall \xi\in\mathbb{R}.
 	$$
 	Note that \(\inf_{x\in\mathbb{R}} z=0\). Then we can choose \(x_0\in\mathbb{R}\) such that
 	$$
 	z(x_0) < \min\left\{\frac{1}{4C},\ \frac{\varepsilon_0}{4\sup_{x\in\mathbb{R}}|a(x)|}\right\}.
 	$$
 	There exists \(R>0\) large enough such that
 	$$
 	\frac{\theta(x-x_0)}{C} \ge z(x_0),\quad \text{for any } x\in B_R^c(0),
 	$$
 	since \(\lim_{|x|\to\infty}\theta(x)=1\).
 	 	Let
 	$$
 	g(x)=z(x)+\frac{\theta(x-x_0)}{C}.
 	$$
 	Then \(g(x) \ge \frac{\theta(x-x_0)}{C} > z(x_0)=g(x_0)\) for \(x\in B_R^c(0)\). Hence \(g(x)\) reaches its minimum at some point \(\overline{x}\in B_R(0)\).
 	Moreover, \(g(\overline{x})\le g(x_0)=z(x_0)\). 	
 	Take
 	$$
 	\phi(x)=z(\overline{x})+\frac{\theta(\overline{x}-x_0)}{C}-\frac{\theta(x-x_0)}{C}.
 	$$
 	We have \(z(x)-\phi(x)=g(x)-g(\overline{x})\ge 0\), and \(z(\overline{x})=\phi(\overline{x})\), hence \(\phi\) touches \(z\) from below at \(\overline{x}\) over \(\mathbb{R}\). Taking \(\phi\) as a test function, we have
 	\begin{align*}
 		-\varepsilon_0
 		&\ge  -\int_{\mathbb{R}}\frac{\theta(\overline{x}-x_0)}{C}J(\overline{x},y)\,dy + a(\overline{x})z(\overline{x}) - q(\overline{x})\frac{\theta'(\overline{x}-x_0)}{C} \\
 		&= -\frac{\mathcal{L}(\theta(\overline{x}-x_0))}{C} + a(\overline{x})g(\overline{x}) \\
 		&\ge -\frac{\varepsilon_0}{4} - \frac{\varepsilon_0}{4} = -\frac{\varepsilon_0}{2}.
 	\end{align*}
 	This is a contradiction. Therefore \(\inf_{\mathbb{R}} z > 0\).
 \end{proof}
 
 \begin{definition}\label{def:convolution}
 	The sup-convolution approximation of an USC function \(u\) is given by
 	$$
 	u^\theta(x) := \sup_{y}\left\{u(y)-\frac{|x-y|^2}{\theta}\right\} = \sup_{y}\left\{u(x+y)-\frac{|y|^2}{\theta}\right\}.
 	$$
 	On the other hand, the inf-convolution approximation of an LSC function \(v\) is given by
 	$$
 	v_\theta(x) := \inf_{y}\left\{v(y)+\frac{|x-y|^2}{\theta}\right\} = \inf_{y}\left\{v(x+y)+\frac{|y|^2}{\theta}\right\}.
 	$$
 \end{definition}
 
 \begin{proposition}(\cite[Proposition 5.4]{CaffarelliSilvestre2009})\label{prop2.3}
 	If \(u\) is bounded and LSC (USC) in \(\mathbb{R}\), then \(u_\theta\) (\(-u^\theta\)) converges to \(u\) (\(-u\)) in the following sense: for any subsequence \(\theta_k\),
 	\begin{enumerate}
 		\item[(a)] for every sequence \(x_k\to x\in\mathbb{R}\), \(\liminf_{k\to\infty} u_{\theta_k}(x_k) \ge u(x)\);
 		\item[(b)] for every \(x\in\mathbb{R}\), there exist \(x_k\to x\) such that \(\limsup_{k\to\infty} u_{\theta_k}(x_k) = u(x)\).
 	\end{enumerate}
 \end{proposition}
 
\begin{proposition}\label{prop:convolution}
	Let \(\phi:\mathbb R\to\mathbb R\) be bounded and uniformly continuous, with modulus of continuity
	$$
	\omega_\phi(\delta)=\sup_{\substack{x,y\in\mathbb R\\ |x-y|\le\delta}}|\phi(x)-\phi(y)|,\qquad \delta\ge 0,
	$$
	\(\phi^\theta\) and $\phi_\theta$ are  the sup- and inf-convolution of $\phi$, respectively. 
	Set
	$$
	r_\theta:=\sup\left\{r\ge0:\ \frac{r^2}{\theta}\le \omega_\phi(r)\right\}.
	$$
	Then  the following hold.\\
	(i) Both \(\phi_\theta(x)\) and \(\phi^\theta(x)\) converge uniformly to \(\phi\). More precisely, 
	$$
	0\le \phi(x)-\phi_\theta(x)\le \omega_\phi(r_\theta)\to0, \text{ as } \theta\to0,
	$$
	$$
	0\le \phi^\theta(x)-\phi(x)\le \omega_\phi(r_\theta)\to0,  \text{ as } \theta\to0.
	$$
	(ii)
		For every \(x\in\mathbb R\), the left and right derivatives of \(\phi_\theta\) and \(\phi^\theta\) exist, and for any $x\in\mathbb R$,
		$$
		|(\phi_\theta)'_\pm(x)|\le \frac{2r_\theta}{\theta},\qquad
		|(\phi^\theta)'_\pm(x)|\le \frac{2r_\theta}{\theta}.
		$$
\end{proposition}

\begin{proof}
	Since
	$$
	\phi^\theta(x)=\sup_{y}\left\{\phi(y)-\frac{|x-y|^2}{\theta}\right\}
	=-\inf_{y}\left\{-\phi(y)+\frac{|x-y|^2}{\theta}\right\}
	=-(-\phi)_\theta(x),
	$$
	and \(-\phi\) is also bounded and uniformly continuous with the same modulus, it suffices to prove the assertions for \(\phi_\theta\). The results for \(\phi^\theta\) follow by applying those for \(\phi_\theta\) to \(-\phi\).
	
	First of all, note that $\omega_\phi$ is continuous since $\phi$ is uniformly continuous. Hence 
	$$
	\frac{r_\theta^2}{\theta}=\omega_\phi(r_\theta), \forall \theta>0.
	$$
	\iffalse
	Let
	$$
	E:=\left\{r\ge0:\ \frac{r^2}{\theta}\le \omega_\phi(r)\right\}.
	$$
	Since \(0\in E\), \(E\neq\emptyset\). Because \(\phi\) is bounded, \(\omega_\phi\) is bounded, say \(\omega_\phi(r)\le M\). If \(r>\sqrt{\theta M}\), then
	$$
	\frac{r^2}{\theta}>M\ge \omega_\phi(r),
	$$
	so \(r\notin E\). Hence \(E\) is bounded above and \(r_\theta=\sup E\) is finite. Also \(r_\theta\ge0\).
	
\(0\le r_\theta<\infty\) and
$$
\frac{r_\theta^2}{\theta}=\omega_\phi(r_\theta).
$$	We claim that
	$$
	\frac{r_\theta^2}{\theta}=\omega_\phi(r_\theta).
	$$
	Indeed, define \(g(r)=r^2/\theta-\omega_\phi(r)\). Since \(\omega_\phi\) and \(r\mapsto r^2/\theta\) are continuous, \(g\) is continuous. If \(g(r_\theta)<0\), then by continuity there exists \(r>r_\theta\) with \(g(r)<0\), i.e.\ \(r\in E\), contradicting the definition of \(r_\theta\). If \(g(r_\theta)>0\), then by continuity there exists \(r<r_\theta\) close to \(r_\theta\) with \(g(r)>0\), so no point in \((r,r_\theta)\) belongs to \(E\), again contradicting the definition of \(r_\theta\). Thus \(g(r_\theta)=0\).
	\fi
	\medskip
It is easy to see that $\phi_\theta(x)=\inf\limits_{y}\left\{\phi(y)+\frac{|x-y|^2}{\theta}\right\} \le \phi(x)$ by taking \(y=x\). 
Choose a minimizer \(y_x\), i.e.,
$$
\phi_\theta(x)=\phi(y_x)+\frac{|x-y_x|^2}{\theta}.
$$
	Suppose \(|x-y|>r_\theta\), then
	\begin{equation}\label{eq:convergece of convolution}
	\phi(x)\geq\phi_\theta(x)=\phi(y_\theta)+\frac{|x-y_\theta|^2}{\theta}
	\ge \phi(x)-\omega_\phi(|x-y_\theta|)+\frac{|x-y_\theta|^2}{\theta}.
	\end{equation}
	By the definition of \(r_\theta\), for \(r>r_\theta\) we have \(r^2/\theta>\omega_\phi(r)\), so the right-hand side is \(>\phi(x)\),
	which is a contradiction. Therefore any minimizing sequence \(y_x\) satisfies
	$$
	|x-y_x|\le r_\theta.
	$$
Moreover, the right hand side of \eqref{eq:convergece of convolution}  great than $\phi(x)-\omega_\phi(|x-y_\theta|)$, which yields that
	$$
	0\le \phi(x)-\phi_\theta(x)\le\omega_\phi(|x-y_\theta|)\le \omega_\phi(r_\theta)=\frac{r_\theta^2}{\theta}.
	$$
	As \(\theta\to0^+\), from \(r_\theta^2/\theta=\omega_\phi(r_\theta)\to0\) we get \(r_\theta\to0\). Since \(\phi\) is uniformly continuous, \(\omega_\phi(r_\theta)\to0\), hence \(\|\phi_\theta-\phi\|_\infty\to0\). The argument for \(\phi^\theta\) is analogous.
	
	Next we estimate the one-sided derivative.
	Note that the function
	$$
	\phi_\theta(x)+\frac{|x|^2}{\theta}
	$$
	is convex, so the left and right derivatives of \(\phi_\theta\) exist everywhere.
	For \(h>0\), by the definition of \(\phi_\theta(x+h)\),
	$$
	\phi_\theta(x+h)\le \phi(y_x)+\frac{|x+h-y_x|^2}{\theta}.
	$$
	Hence
	$$
	\phi_\theta(x+h)-\phi_\theta(x)
	\le \frac{1}{\theta}\left(|x+h-y_x|^2-|x-y_x|^2\right)
	=\frac{2h}{\theta}(x-y_x)+\frac{h^2}{\theta}.
	$$
	Dividing by \(h>0\) and letting \(h\to0^+\),
	$$
	(\phi_\theta)'_+(x)\le \frac{2(x-y_x)}{\theta}\le \frac{2r_\theta}{\theta}.
	$$
	
	On the other hand, take a minimizer \(y_{x+h}\) for \(x+h\). By the definition of \(\phi_\theta(x)\),
	$$
	\phi_\theta(x)\le \phi(y_{x+h})+\frac{|x-y_{x+h}|^2}{\theta}.
	$$
We get
	$$
	\phi_\theta(x+h)-\phi_\theta(x)
	\ge \frac{2h}{\theta}(x-y_{x+h})+\frac{h^2}{\theta}=\frac{2h}{\theta}(x+h-y_{x+h})-\frac{h^2}{\theta}.
	$$
	Dividing by \(h>0\), letting \(h\to0^+\), and using \(|x+h-y_{x+h}|\le r_\theta\), we have
	$$
	(\phi_\theta)'_+(x)\ge -\frac{2r_\theta}{\theta}.
	$$
	To sum up
	$$
	|(\phi_\theta)'_+(x)|\le \frac{2r_\theta}{\theta}.
	$$
	Analogous arguments yield the corresponding results for \(\phi^\theta\).
\end{proof}

 \begin{proof}[Proof of Theorem \ref{thm:main_result1}]
 	Since \(\mathcal{A}_1\subset\mathcal{A}_2\), then (i) follows from the definition of the generalized principle eigenvalues.

 	Let us now show that  (ii) holds.
 	We prove only the case $\nu=0$; the case $\nu\neq 0$ can be proved similarly.
 	A simple argument establishes that \(\underline{\lambda}(0,\mathcal{A}_1) \le \overline{\lambda}(0,\mathcal{A}_1)\). Indeed, if not, then  \(\underline{\lambda}(0,\mathcal{A}_1) > \overline{\lambda}(0,\mathcal{A}_1)\), which yields that there exist \(\lambda, \delta_0 > 0\) such that \(\underline{\lambda}(0,\mathcal{A}_1) > \lambda >\lambda-2\delta_0 > \overline{\lambda}(0,\mathcal{A}_1)\), and there exist \(\psi, \phi \in \mathcal{A}_1\) satisfying the following in the classical sense: 	
 	\begin{equation}\label{eq:sub and super}
 		\begin{cases}
 			\mathcal{L}\psi \ge \lambda\psi, \\
 			\mathcal{L}\phi \le (\lambda-2\delta_0)\phi,
 		\end{cases}
 	\end{equation}
 	Therefore, setting \(\gamma = \inf\frac{\phi}{\psi} > 0\), we obtain \(\mathcal{L}(\phi-\gamma\psi) \le \lambda(\phi-\gamma\psi)-2\delta_0\phi\). This contradicts Lemma \ref{lem2.1}.
 	
 	Let us now prove   \(\underline{\lambda}(0,\mathcal{A}_2) \le \overline{\lambda}(0,\mathcal{A}_2)\) by contradiction.
 	Suppose  there exist \(\lambda,\delta_0>0\) such that \(\underline{\lambda}(0,\mathcal{A}_2)>\lambda>\lambda-2\delta_0>\overline{\lambda}(0,\mathcal{A}_2)\) and \(\psi,\phi\in\mathcal{A}_2\) such that \eqref{eq:sub and super} holds	in the viscosity sense.
 	To apply Lemma \ref{lem2.1}, we first state the following: %Define \(\gamma=\inf_{\mathbb{R}}\frac{\phi}{\psi}>0\) and \(z=\phi-\gamma\psi\ge 0\), with \(\inf_{\mathbb{R}} z=0\).
 	\\ \textbf{Claim:} It holds that
 	$$
 	\begin{cases}
 		\mathcal{L}(\psi^\theta) \ge \lambda\psi^\theta -  c_\theta^1,\\
 		\mathcal{L}(\phi_\theta) \le (\lambda-2\delta_0)\phi_\theta + c_\theta^2,
 	\end{cases}
 	$$
 	where \(c_\theta^1,c_\theta^2\to 0\) as \(\theta\to 0\), depending on the moduli of continuity of \(\psi\) and \(\phi\), respectively, in the viscosity sense.\\	
	\textbf{Proof of Claim:} We only prove the first inequality and the proof of the second one is similar. 
	Suppose \(\Phi\) touches \(\psi^\theta\) from above at \(x_0\). Then
	\begin{equation}
		\Phi(x) \ge \psi^\theta(x) \quad \text{for } x \text{ in some small neighborhood }  \text{ of } x_0. \tag{24}
	\end{equation}
	Let \(y_{x_0}\) be a minimizer, i.e.,
	$$
	\Phi(x_0) = \psi^\theta(x_0) = \psi(y_{x_0}) - \frac{(x_0 - y_{x_0})^2}{\theta},
	$$
	and set \(\Psi(x) = \Phi(x + x_0 - y_{x_0}) + \frac{(x_0 - y_{x_0})^2}{\theta}\). Then \(\Phi'(x_0) = \Psi'(y_{x_0})\),
	\begin{align*}
		\Psi(y_{x_0}) &= \Phi(y_{x_0} + x_0 - y_{x_0}) + \frac{(x_0 - y_{x_0})^2}{\theta} \\
		&= \Phi(x_0) + \frac{(x_0 - y_{x_0})^2}{\theta} = \psi(y_{x_0})
	\end{align*}
		Furthermore, in a neighborhood of \(y_{x_0}\), we have
	\begin{align*}
		\Psi(x) &= \Phi(x + x_0 - y_{x_0}) + \frac{(x_0 - y_{x_0})^2}{\theta} \ge \psi^\theta(x + x_0 - y_{x_0}) + \frac{(x_0 - y_{x_0})^2}{\theta} \\
		&= \sup_y \left\{ \psi(y) - \frac{(x + x_0 - y_{x_0} - y)^2}{\theta} \right\} + \frac{(x_0 - y_{x_0})^2}{\theta} \\
		&\ge \psi(x) - \frac{(x_0 - y_{x_0})^2}{\theta} + \frac{(x_0 - y_{x_0})^2}{\theta} = \psi(x)
	\end{align*}
	where the first inequality uses equation (24). The above shows that \(\Psi\) touches \(\psi\) from above at \(y_{x_0}\). Therefore, from the equation satisfied by \(\psi\), we have:
	$$
	\int_{\mathbb{R}} \psi(y) J(y_{x_0}, y) \, dy + (a(y_{x_0}) - \lambda) \psi(y_{x_0}) + q(y_{x_0}) \Psi'(y_{x_0}) \ge 0.
	$$
	Thus
	\begin{align*}
		&\quad\int_{\mathbb{R}} \psi^\theta(y) J(x_0, y) \, dy + (a(x_0) - \lambda) \psi^\theta(x_0) + q(x_0) \Phi'(x_0) \\
		&\ge \int_{\mathbb{R}} \psi(y) J(x_0, y) \, dy + (a(x_0) - \lambda) \left( \psi(y_{x_0}) - \frac{(x_0 - y_{x_0})^2}{\theta} \right) + q(x_0) \Phi'(x_0) \\
		&= \int_{\mathbb{R}} \psi(y) J(y_{x_0}, y) \, dy + (a(y_{x_0}) - \lambda) \psi(y_{x_0}) + q(y_{x_0}) \Psi'(y_{x_0}) \\
		&\quad + \int_{\mathbb{R}} \psi(y) (J(x_0, y) - J(y_{x_0}, y)) \, dy - (a(x_0) - \lambda) \frac{(x_0 - y_{x_0})^2}{\theta} \\
		&\quad + (a(x_0) - a(y_{x_0})) \psi(y_{x_0}) + (q(x_0) - q(y_{x_0})) \Phi'(x_0) \\
		&\ge \int_{\mathbb{R}} \psi(y) (J(x_0, y) - J(y_{x_0}, y)) \, dy - \|a - \lambda\|_\infty \frac{(x_0 - y_{x_0})^2}{\theta} \\
		&\quad + \|\psi\|_\infty \omega_a(|x_0 - y_{x_0}|) - L_q |x_0 - y_{x_0}| \cdot |\Phi'(x_0)|,
	\end{align*}	
	where $L_q$ is the Lipschitz constant of $q$.
	Note that \(\Phi\) and \(\psi^\theta\) touch at \(x_0\). Thus \(\Phi'(x_0)\) belongs to the subdifferential of \(\psi^\theta\). Therefore, combined with (ii) of Proposition 3.4, we have $	|\Phi'(x_0)| \le \frac{2|x_0 - y_{x_0}|}{\theta},$
	which implies
	$$
	|x_0 - y_{x_0}| \cdot |\Phi'(x_0)| \le \frac{2(x_0 - y_{x_0})^2}{\theta}.
	$$
	Besides, 
		\begin{align*}
			&\quad \left|\int_{\mathbb{R}} \psi(y) (J(x_0, y) - J(y_{x_0}, y)) \, dy\right| \\
		&\le \left|\int_{\mathbb{R}} (\psi(x_0-y)-\psi(y_{x_0}-y)) J(x_0, x_0-y)\, dy\right|\\
		&\quad +\left|\int_{\mathbb{R}} \psi(z_0-y) (J(x_0,x_0-y) - J(y_{x_0},y_{x_0}-y)) \, dy\right|  \\
		&\le \int J_0(y)dy\omega_{\psi}(|x_0-y_{x_0}|)+o(|x_0-y_{x_0}|),
	\end{align*}	
	where (J2) was used in the last inequality.
	Hence, the right-hand side of the inequality above is \(\ge -c_\theta^1\), where \(c_\theta^1\) is independent of \(x_0\) and satisfies
	$$
	c_\theta^1 \to 0 \quad \text{as} \quad \theta \to 0.
	$$
	Consequently,
	$$
	\mathcal{L}\psi^\theta \ge \lambda \psi^\theta - c_\theta^1.
	$$
	This completes the proof of claim.
	
	Noting that \(\phi_\theta \rightrightarrows \phi\) and \(\psi^\theta \rightrightarrows \psi\), we can choose \(\theta\) sufficiently small
	such that
	$$
	\gamma c_\theta^1 + c_\theta^2 - 2\delta_0 \inf \phi_\theta < 0,
	$$
	where the constant is defined as \(\gamma := \inf\limits_x \frac{\phi_\theta}{\psi^\theta}\).
 We will show that
 	$$
 	\mathcal{L}(\phi_\theta - \gamma\psi^\theta) \le \lambda(\phi_\theta - \gamma\psi^\theta) + \gamma c_\theta^1 + c_\theta^2 - 2\delta_0 \phi_\theta
 	$$
 	in the viscosity sense.
 	
 	Suppose \(\varphi \in C^1\) touches \(\phi_\theta - \gamma\psi^\theta\) from below at \(x_0\). Then the derivatives of \(\phi_\theta\) and \(\psi^\theta\) exist at \(x_0\). In fact, since \(\phi_\theta\) and \(-\gamma\psi^\theta\) are semiconcave, there exists \(M\) large enough such that in a neighborhood \(U\) of \(x_0\), the functions
 	$$
 	F(x) := \phi_\theta(x) - \frac{M}{2}(x-x_0)^2, \quad G(x) := -\gamma\psi^\theta - \frac{M}{2}(x-x_0)^2
 	$$
 	are concave, and \(H(x) := \varphi(x) - M(x-x_0)^2\) touches \(F(x)+G(x)\) from below at \(x_0\).
 	
 	 	We only need to prove that the derivatives of \(F\) and \(G\) exist at \(x_0\). Let \(\partial F(x_0)\) and \(\partial G(x_0)\) denote the subdifferentials of \(F\) and \(G\) at \(x_0\), respectively. Then in a neighborhood \(U\) of \(x_0\),
 	$$
 	\begin{cases} 
 		F(x) \le F(x_0) + p(x-x_0) \\ 
 		G(x) \le G(x_0) + q(x-x_0) 
 	\end{cases}
 	$$
 	where \(p \in \partial F(x_0)\), \(q \in \partial G(x_0)\). Thus,
 	$$
 	F(x) + G(x) \le F(x_0) + G(x_0) + (p+q)(x-x_0).
 	$$
 	On the other hand, we have \(F(x_0) + G(x_0) = H(x_0)\) and
 	$$
 	F(x) + G(x) \ge H(x) = H(x_0) + H'(x_0)(x-x_0) + o(|x-x_0|).
 	$$
 	Therefore,
 	$$
 	H'(x_0)(x-x_0) + o(|x-x_0|) \le (p+q)(x-x_0),
 	$$
 	that is, \((H'(x_0) - p - q)(x-x_0) \le o(|x-x_0|)\),
 	which yields that
 	$$
 	H'(x_0) = p + q.
 	$$
 	This shows that both \(\partial F(x_0)\) and \(\partial G(x_0)\) are singletons. Hence, \(F\) and \(G\) are differentiable at \(x_0\). Consequently, \(\phi_\theta\) and \(\psi^\theta\) are differentiable at \(x_0\).
 	
 	Since \(\phi_\theta\) and \(\psi^\theta\) are differentiable at \(x_0\), we have
 	$$
 	\begin{cases} 
 		\mathcal{L}(\gamma\psi^\theta) \ge \lambda\gamma\psi^\theta - \gamma c_\theta^1 \\ 
 		\mathcal{L}\phi_\theta \le (\lambda - 2\delta_0)\phi_\theta + c_\theta^2 
 	\end{cases}
 	$$
 	holding at \(x_0\) in the classical sense.
 	Thus,
 	$$
 	\mathcal{L}(\phi_\theta - \gamma\psi^\theta) \le \lambda(\phi_\theta - \gamma\psi^\theta) + \gamma c_\theta^1 + c_\theta^2 + 2\delta_0 \phi_\theta.
 	$$
 	 	Since \(\inf(\phi_\theta - \gamma\psi^\theta) = 0\), the above inequality contradicts Lemma \ref{lem2.1}.
 	 	Therefore, \(\underline{\lambda}(0, \mathcal{A}_2) \le \overline{\lambda}(0, \mathcal{A}_2)\).
 \end{proof}

\subsection{Proof of Theorem \ref{thm:main_result2}}

\begin{proof}[Proof of Theorem \ref{thm:main_result2}]	
	Proof of (i). 
	We first prove that the integral operator $(T\phi)(x)=\int_{\Omega}J(x,y)\phi(y)\,dy$
	is compact from $C^1(\Omega)\cap C(\overline{\Omega})$ to $C(\overline{\Omega}).$
	Take any bounded set \(B\subset C(\overline{\Omega})\), i.e., there exists a
	constant \(M_B>0\) such that
	$$
	\|\phi\|_\infty\le M_B,\quad\forall\phi\in B.
	$$
%	To prove that \(T\) is compact, it suffices to show that the image set 	\(T(B)=\{T\phi:\phi\in B\}\) is uniformly bounded and equicontinuous in \(C(\overline{\Omega})\).	
	By the upper bound in (J1), it is easy to see that \(T(B)\) is uniformly bounded.	
	Next, we prove equicontinuity. For any \(\phi\in B\) and
	\(x_1,x_2\in\overline{\Omega}\), by the linearity of \(T\),
	\begin{equation}\label{eq:equi}
		|(T\phi)(x_1)-(T\phi)(x_2)|
		\le M_B\int_{\Omega}\bigl|J(x_1,y)-J(x_2,y)\bigr|\,dy.
	\end{equation}
	By the lower bound in (J1), \(J(x_i,y)\) can be nonzero only when
	\(|x_i-y|<\delta_0\), so the integral is actually taken over
	$$
	D=\{y\in\Omega:|x_1-y|<\delta_0\text{ or }|x_2-y|<\delta_0\}.
	$$
	Split \(D\) into the overlapping part
	\(D_{\mathrm{over}}=\{y:|x_1-y|<\delta_0\text{ and }|x_2-y|<\delta_0\}\)
	and the non-overlapping part
	\(D_{\mathrm{diff}}=D\setminus D_{\mathrm{over}}\).	
	On \(D_{\mathrm{over}}\), set \(y=x_1-z\). Then
	$$
	J(x_1,y)=J(x_1,x_1-z),\qquad
	J(x_2,y)=J(x_2,x_2-(z+(x_2-x_1))).
	$$
	The difference of the two translation parameters is exactly \(x_2-x_1\).
	By (J2), for any \(\varepsilon>0\), there exists \(\delta>0\) such that
	whenever \(|x_1-x_2|<\delta\),
	$$
	|J(x_1,y)-J(x_2,y)|<\frac{\varepsilon}{2M_B|\Omega|},
	$$
	so the integral over the overlapping region is less than
	\(\varepsilon/2\).	
	On \(D_{\mathrm{diff}}\), its measure is at most \(C|x_1-x_2|\)
	(the constant \(C\) depends only on \(\delta_0\)), and by the upper bound
	in (J1), \(J\) is bounded on compact sets; denote this bound by \(M_0\).
	Hence
	$$
	\int_{D_{\mathrm{diff}}}|J(x_1,y)-J(x_2,y)|\,dy
	\le 2M_0\cdot C|x_1-x_2|.
	$$
	Taking \(|x_1-x_2|\) sufficiently small makes this integral less than
	\(\varepsilon/2\) as well.
	Combining the two parts, whenever \(|x_1-x_2|<\delta\),
	$$
	|(T\phi)(x_1)-(T\phi)(x_2)|<\varepsilon,
	$$
	and this estimate holds uniformly for all \(\phi\in B\). Therefore
	\(T(B)\) is equicontinuous.	
    To sum up, \(T(B)\) is uniformly bounded and equicontinuous. 
    By the Arzel\`a--Ascoli theorem,
	% \(T(B)\) is relatively compact in \(C(\overline{\Omega})\). Since \(B\) is an arbitrary bounded set and	\(T\) maps bounded sets to relatively compact sets, 
	\(T\) is a compact 	operator.	    
     Using this compactness, we can establish results analogous to those in \cite[Theorem 2.1]{Coville2020} (see also \cite[Theorem 4.1]{LCW}) and \cite[Proposition 1.1]{Coville2020}, and then derive (i) by following the argument in \cite[Theorem 1.2]{Coville2020}.
	
	Proof of (ii). With (i) at hand, and using the Harnack inequality established in this paper together with the Remark \ref{re:more-about-HIE}, (ii) can be proved in a manner analogous to  \cite[Theorem 1.4]{Coville2020}.
	
	Proof of (iii).
	Fix arbitrary $M>0$. Choose $\mu \in \mathbb{R}$ with $|\mu|$ sufficiently large such that
$$
\mu q(x) + a(x) \ge M + 1,\quad \forall x \in \Omega.
$$
This is possible because $\inf_\Omega |q| > 0$ by assumption (A2). If $q>0$ on $\Omega$, take $\mu>0$; if $q<0$, take $\mu<0$.
Let $\phi(x) = e^{\mu x}\in \mathcal A_\Omega$. 
Then for $x \in \Omega$,
$$
\mathcal L\phi(x) = \int_\Omega e^{\mu y}J(x,y)\,dy + a(x)e^{\mu x} + q(x)\mu e^{\mu x} \ge (a(x)+\mu q(x))e^{\mu x} = (M+1)\phi(x).
$$
Thus $\phi$ satisfies $\mathcal L\phi \ge (M+1)\phi$ in $\Omega$, which yields that
$\underline{\lambda}(0,\mathcal A_\Omega) \ge M+1.$
Since $M>0$ was arbitrary,
$$
\underline{\lambda}(0,\mathcal A_\Omega) = +\infty.
$$
On the other hand, it is easy to see that  $\overline{\lambda}(0,\mathcal A_\Omega) < +\infty$ by taking $\phi \equiv 1\in \mathcal A_\Omega$ as a test function.
Therefore,
$$
\overline{\lambda}(0,\mathcal A_\Omega) \le C < +\infty = \underline{\lambda}(0,\mathcal A_\Omega).
$$
This complete the proof.
\end{proof}

\subsection{Proof of Theorem \ref{the:combining hamel and du}}
\begin{proof}[Proof of Theorem \ref{the:combining hamel and du}]
	If $\mathcal I=\emptyset$, then there is nothing need to prove since all the quantities in \eqref{eq:1 of the:combining hamel and du} are $+\infty$. Hence, from now on, we assume that $\mathcal I\ne\emptyset$.
	
		It is easy to see that $\overline{\lambda}(0,\mathcal{A}_{\mathbb{R}}) \le \inf\limits_{\nu\in\mathbb{R}}\overline{\lambda}(\nu,\mathcal{A}_1)$.
	In fact, for any $\nu\in\mathbb R$ and $\varepsilon>0$, there exists $\phi\in\mathcal A_1$ such that $\mathcal L_\nu\phi\leq(\overline{\lambda}(\nu,\mathcal{A}_1)+\varepsilon)\phi$.
	Note that $\phi\in\mathcal A_1\subset \mathcal A_{\R}$.
	The definition of $\overline{\lambda}(0,\mathcal{A}_{\mathbb{R}})$ yields that 
	$$\overline{\lambda}(0,\mathcal{A}_{\mathbb{R}}) \le \overline{\lambda}(\nu,\mathcal{A}_1)+\varepsilon, \forall \nu\in\mathbb R, \varepsilon>0.$$
	Hence $\overline{\lambda}(0,\mathcal{A}_{\mathbb{R}}) \le \inf\limits_{\nu\in\mathbb{R}}\overline{\lambda}(\nu,\mathcal{A}_1)$.
	Besides, it is easy to see from the definition of  $\overline{\lambda}(0,\mathcal{A}_{\Omega})$ that $\overline{\lambda}(0,\mathcal{A}_{\Omega})\le \overline{\lambda}(0,\mathcal{A}_{\mathbb{R}})$ for any $\Omega\subset\mathbb R,$
	which yields that $\lim\limits_{n\to\infty}\overline{\lambda}(0,\mathcal{A}_{\Omega_n})\leq\overline\lambda(0,\mathcal A_\mathbb R).$
	
	To prove $\inf\limits_{\nu\in\mathbb{R}}\underline{\lambda}(\nu,\mathcal{A}_1)\le\lim\limits_{n\to\infty}\overline{\lambda}(0,\mathcal{A}_{\Omega_n})$, we first show that $ \inf\limits_{\nu\in\mathbb{R}}\underline{\lambda}(\nu,\mathcal{A}_1)\le \overline{\lambda}(0,\mathcal{A}_{\mathbb{R}})$ by  homogenization.
	There are two cases we need to consider:
	
	Case 1: \(J\) satisfies (J1).
	Under this assumption, by (ii) of Theorem \ref{thm:main_result2}, there exists $\Psi \in \mathcal A_{\mathbb R}$ such that
	$$
	\int_{\mathbb R} J(x,y)\Psi(y)\,dy + q(x)\Psi(x) + a(x)\Psi(x) = \lambda \Psi(x),
	$$
	where
	$$
	\lambda=\overline{\lambda}(0,\mathcal A_{\mathbb R}) = \inf\{\lambda \mid \exists \Psi \in \mathcal A_{\mathbb R} \text{ such that } \mathcal L_\nu\Psi \ge \lambda \Psi \text{ in } \mathbb R\}.
	$$
	We normalize $\Psi(0)=1$.	Let
	$$
	w_\varepsilon(x) = \varepsilon \ln \Psi\left(\frac{x}{\varepsilon}\right), \qquad w_\varepsilon(0)=0.
	$$
	
	\textbf{Claim 1:} $w_\varepsilon(x)$ is locally uniformly bounded, and $\{w'_\varepsilon\}$ is uniformly bounded in $\varepsilon$.
	
	Without loss of generality, assume $x<y$. By Theorem \ref{harnack}, the Harnack type inequality, there exists $C>0$ such that
	$$
	\frac{\Psi(\frac{x}{\varepsilon})}{\Psi(\frac{x}{\varepsilon}+1)}
	\cdot \frac{\Psi(\frac{x}{\varepsilon}+1)}{\Psi(\frac{x}{\varepsilon}+2)}
	\cdots
	\frac{\Psi(\frac{x}{\varepsilon}+\frac{y}{\varepsilon})}{\Psi(\frac{y}{\varepsilon})}
	\le C^{\left( \left| \frac{y-x}{\varepsilon} \right| + 1 \right)}.
	$$
	Consequently,
	\begin{equation}\label{eq:linear growth}
		|w_{\varepsilon}(x)-w_{\varepsilon}(y)|
		=\left|\varepsilon \ln \Psi(\frac{x}{\varepsilon}) - \varepsilon \ln \Psi(\frac{y}{\varepsilon})\right|
		\le \varepsilon C\left( \left| \frac{y-x}{\varepsilon} \right| + 1 \right).
	\end{equation}
	Hence $w_\varepsilon$ is locally uniformly bounded.
	
	The function $w_\varepsilon(x)$ satisfies
	$$
	\int_{\mathbb{R}} J\left(\frac{x}{\varepsilon}, y\right) e^{\frac{w_\varepsilon(\varepsilon y)}{\varepsilon}} \, dy
	+ q\left(\frac{x}{\varepsilon}\right) e^{\frac{w_\varepsilon(x)}{\varepsilon}} w'_\varepsilon(x)
	+ a\left(\frac{x}{\varepsilon}\right) e^{\frac{w_\varepsilon(x)}{\varepsilon}}
	= \lambda e^{\frac{w_\varepsilon(x)}{\varepsilon}}
	$$
	i.e.
	\begin{equation}\label{eq:lambda}
		\int_{\mathbb{R}} J\left(\frac{x}{\varepsilon}, y\right) e^{\frac{w_\varepsilon(\varepsilon y) - w_\varepsilon(x)}{\varepsilon}} \, dy
		+ q\left(\frac{x}{\varepsilon}\right) w^\prime_\varepsilon(x)
		+ a\left(\frac{x}{\varepsilon}\right)
		= \lambda,
	\end{equation}
	which yields that 
	\begin{equation}\label{eq:wprime}
		w^\prime_\varepsilon(x)=\frac{\lambda- a\left(\frac{x}{\varepsilon}\right)}{q\left(\frac{x}{\varepsilon}\right)}-\frac{1}{q\left(\frac{x}{\varepsilon}\right)}\int_{\mathbb{R}} J\left(\frac{x}{\varepsilon}, y\right) \exp\{\frac{w_\varepsilon(\varepsilon y) - w_\varepsilon(x)}{\varepsilon}\} dy.
	\end{equation}
	By \eqref{eq:linear growth}, we have
	$$
	\int_{\mathbb R} J(\frac{x}{\varepsilon},y) \exp\left(\frac{w(\varepsilon y)-w(x)}{\varepsilon}\right)dy
	\le \int_{\mathbb R} J(\frac{x}{\varepsilon},y) e^{c(|\frac{x}{\varepsilon}-y|)+c}dy < \infty.
	$$
	Thus $w^\prime_{\varepsilon}$  is bounded by \eqref{eq:wprime}, and hence (along a subsequence)
	$$
	w_{\varepsilon_n}(x) \to w(x) \quad \text{locally uniformly}.
	$$
	Moreover, $w$ is Lipschitz on $\mathbb R$ and grows at most linearly.	
	Take a smooth function $\chi$ such that as $x\to 0$, $\chi(x)$ tends to $-\infty$ superlinearly. Suppose
	$$
	w(x)-\chi(x)
	$$
	attains a global strict minimum at $x_0$ (strictness is achievable). Let $\nu_0=\chi(x_0)$, and 
	\iffalse
	suppose
	$$
	w_{\varepsilon_n}(x)-\chi(x)
	$$
	attains a global minimum at $\xi_n$.
	Then there exists a subsequence (still denoted by $\xi_n$) $\xi_n$ such that as $n\to\infty$,
	$$
	J(x+\xi_n, y+\xi_n) \to J^*(x,y),
	$$
	$$
	q(x+\xi_n) \to q^*(x), \qquad a(x+\xi_n)\to a^*(x),
	$$
	since $q, a$ are almost periodic, and  the family of functions $\{J(\cdot, \cdot + y) : y \in \mathbb R\}$ is uniformly almost periodic in $y$, i.e.
	$$
	\forall \eta > 0,\ \exists \text{ relatively dense } T \subset \mathbb R,\ \forall t \in T,\ \forall x, y \in \mathbb R:
	|J(x+t, x+t+y) - J(x, x+y)| < \eta.
	$$		Moreover, 
	\fi
	for any $\eta>0$, there exists $\phi\in\mathcal A_1$ such that
	\begin{equation}\label{eq:lambda star}
		\mathcal L_{\nu_0}\phi= \int_{\mathbb R} J(\cdot,y)e^{-\nu_0(\cdot-y)}\phi(y)dy +q\phi^\prime  + (a+\nu_0 q)\phi \ge (\underline{\lambda}(\nu_0,\mathcal A_1)-\eta)\phi.
	\end{equation}
	%	where
	%	$$
	%	\mathcal L^*_\nu\phi = \int_{\mathbb R} J^*(x,y)e^{-\nu(x-y)}\phi(y)dy +q^*(x)\phi^\prime (x) + (a^*(x)+\nu q^*(x))\phi(x).
	%	$$
	
	%	\textbf{Claim 2:} One can choose $\phi$ such that $\|\phi'(s)/\phi(s)\|$ is bounded on $\mathbb R$.
	
	%	\textbf{Proof:} Take $\phi=e^{u_\delta}$ with $\delta$ small (depends on $\eta$), where $e^{u_\delta}$ satisfies
	%	$$
	%	\mathcal L^*_\nu e^{u_\delta} = \delta u_\delta e^{u_\delta}.
	%	$$
	%	Then the boundedness of $u^\prime_\delta=\phi'/\phi$ follows analogously to the proof of boundedness of $w'$ in Claim 1.
	
	Now consider
	$$
	w_{\varepsilon}(x) - \left[\chi(x) + \varepsilon \ln \phi\left(\frac{x}{\varepsilon}\right)\right].
	$$
	It attains a global minimum at ${x}_\varepsilon$.
	Then one can find  a subsequence $\{x_n\}$ of $\{x_\varepsilon\}$ such that $\lim\limits_{n\to\infty}x_n= x_0$. Moreover,
	$$w_{\varepsilon_n}(\varepsilon_n y) - \left[\chi(\varepsilon_n y) + \varepsilon_n \ln \phi(y)\right] \ge  w_{\varepsilon_n}( x_n) - \left[\chi(x_n) + \varepsilon_n \ln \phi(\frac{ x_n}{\varepsilon_n})\right].$$
	Hence
	$$w_{\varepsilon_n}(\varepsilon_n y) -w_{\varepsilon_n}(x_n) \ge  \chi(\varepsilon_n y)- \chi(x_n) + \varepsilon_n \ln \frac{\phi(y)}{\phi(\frac{ x_n}{\varepsilon_n})},$$
	and 
	$$
	w'_{\varepsilon_n}(x_n)= \chi'(x_n) + \frac{\phi^\prime(\frac{x_n}{\varepsilon_n})}{\phi(\frac{x_n}{\varepsilon_n})}.
	$$
	
	Substituting into \eqref{eq:lambda} gives
	$$
	\lambda = \int_{\mathbb R} J(\frac{x_n}{\varepsilon_n}, y)
	\exp\left\{\frac{\chi(\varepsilon_n y)-\chi(x_n)}{\varepsilon_n}\right\}
	\frac{\phi(y)}{\phi(\frac{x_n}{\varepsilon_n})}\,dy
	$$
	$$
	+ q(\frac{ x_n}{\varepsilon_n})\left[\chi'(x_n) + \frac{\phi'(\frac{x_n}{\varepsilon_n})}{\phi(\frac{x_n}{\varepsilon_n})}\right]
	+ a(\frac{ x_n}{\varepsilon_n}).
	$$
	Note that $x_n\to x_0$.
	Then
	%	\begin{equation}
		%		\left|a\left(\frac{\hat x_n}{\varepsilon_n}\right) - a^*({s}_n) \right|= \left|a(s_n + \phi_{\varepsilon_n}) - a^*({s}_n)\right| \to 0,
		%	\end{equation}
	\begin{equation}
		\left| q\left(\frac{ x_n}{\varepsilon_n}\right) \chi'(\hat x_n) - q\left(\frac{ x_n}{\varepsilon_n}\right) \nu_0 \right| \to 0,
	\end{equation}
	and for any  \( R > 0 \) fixed, we have
	$$
	\left|\frac{\chi(\varepsilon_n y ) - \chi(x_n )}{\varepsilon_n} - \nu_0(y - \frac{x_n}{\varepsilon_n})\right| \to 0
	$$
	as $n\to\infty$, uniformly in $y\in B_R(\frac{x_n}{\varepsilon_n})$.
	
	Now take $R$  sufficiently large such that
	$$
	\int_{B_R^c(\frac{x_n}{\varepsilon_n})} J(\frac{x_n}{\varepsilon_n},y) e^{-\nu_0(\frac{x_n}{\varepsilon_n}-y)}\frac{\phi(y)}{\phi(\frac{x_n}{\varepsilon_n})}dy < \eta,
	$$
	and for such fixed $R$
	$$
	\begin{aligned}
		& \lim\limits_{n\to\infty}\Bigg|\int_{B_R(\frac{x_n}{\varepsilon_n})} J(\frac{x_n}{\varepsilon_n}, y)\exp\left\{\frac{\chi(\varepsilon_n y)- \chi(x_n)}{\varepsilon}\right\}	\frac{\phi(y)}{\phi(\frac{x_n}{\varepsilon_n})}dy\\
		&\quad\quad- \int_{B_R(\frac{x_n}{\varepsilon_n})} J(\frac{x_n}{\varepsilon_n},y)e^{-\nu_0(\frac{x_n}{\varepsilon_n}-y)}\frac{\phi(y)}{\phi(\frac{x_n}{\varepsilon_n})}dy\Bigg|= 0.
	\end{aligned}
	$$
	Hence, for $n$  sufficiently large,
	\begin{equation}\label{eq:estimate3}
		\begin{aligned}
			& \quad \int_{\mathbb R} J(\frac{x_n}{\varepsilon_n}, y)\exp\left\{\frac{\chi(\varepsilon_n y)-\chi(x_n)}{\varepsilon_n}\right\}
			\frac{\phi(y)}{\phi(\frac{x_n}{\varepsilon_n})}\,dy\\
			& \ge \int_{B_R(\frac{x_n}{\varepsilon_n})} J(\frac{x_n}{\varepsilon_n}, y)\exp\left\{\frac{\chi(\varepsilon_n y)- \chi(x_n)}{\varepsilon}\right\}	\frac{\phi(y)}{\phi(\frac{x_n}{\varepsilon_n})}dy\\
			& \ge \int_{B_R(\frac{x_n}{\varepsilon_n})} J(\frac{x_n}{\varepsilon_n},y)e^{-\nu_0(\frac{x_n}{\varepsilon_n}-y)}\frac{\phi(y)}{\phi(\frac{x_n}{\varepsilon_n})}dy-\eta\\
			& \ge \int_{\mathbb R} J(\frac{x_n}{\varepsilon_n},y)e^{-\nu_0(\frac{x_n}{\varepsilon_n}-y)}\frac{\phi(y)}{\phi(\frac{x_n}{\varepsilon_n})}dy-2\eta.
		\end{aligned}
	\end{equation}	
	Combined with \eqref{eq:lambda star}-\eqref{eq:estimate3}, for $n$ sufficiently large,
	\begin{equation}
		\begin{split}
			\lambda &\ge \int_{\mathbb R} J(\frac{x_n}{\varepsilon_n},y) e^{-\nu_0(\frac{x_n}{\varepsilon_n}-y)} \frac{\phi(y)}{\phi(\frac{x_n}{\varepsilon_n})}dy
			+ q(\frac{x_n}{\varepsilon_n})\left(\nu_0 +\frac{\phi^\prime(\frac{x_n}{\varepsilon_n})}{\phi(\frac{x_n}{\varepsilon_n})}\right)+ a(\frac{x_n}{\varepsilon_n}) - 3\eta\\
			&\ge \underline\lambda(\nu_0,\mathcal A_1)-4\eta
		\end{split}
	\end{equation}
	Since $\eta>0$ is arbitrary, we obtain $	\lambda \ge \underline\lambda(\nu_0,\mathcal A_1)\ge\inf\limits_{\nu\in\R} \underline\lambda(\nu,\mathcal A_1).$
	
	Case 2: \(J\) satisfies (J1)$^\prime$. Let
	$$
	J_k(x,y) = J(x,y)\eta_k(x-y),
	$$
	where $\eta_k(x)=1$ if $x\in([-k,k])$, $\in(0,1)$ if $|x|\in(k,k+1/2)$, $=0$ if $|x|\ge k+1/2$. 
	Denote the generalized principal eigenvalues defined by \(J_k\) as \(\overline{\lambda}_k(\nu, \mathcal{A})\) and \(\underline{\lambda}_k(\nu, \mathcal{A})\).
	
	Since \(J_k \le J\), it is obvious from the definition that \(\overline{\lambda}(0, \mathcal{A}_{\mathbb{R}}) \ge \overline{\lambda}_k(0, \mathcal{A}_{\mathbb{R}})\) holds for any \(k\in\mathbb Z_+\). From the result of the previous step, we have
	\(
	\overline{\lambda}_k(0, \mathcal{A}_{\mathbb{R}}) \ge \inf_{\nu \in \mathbb{R}} \underline{\lambda}_k(\nu, \mathcal{A}_1),\ \forall k\in\mathbb Z_+.
	\)
	Therefore, it suffices to prove that 
	\begin{equation}\label{eq:star}
	\lim\limits_{k\to+\infty}	\inf_{\nu \in \mathbb{R}} \underline{\lambda}_k(\nu, \mathcal{A}_1) = \inf_{\nu \in \mathbb{R}} \underline{\lambda}(\nu, \mathcal{A}_1). 
	\end{equation}
	Since \(\kappa_0 \mathbf{1}_{(x-\delta_0,x+\delta_0)}(y) \leq J(x,y)\), then
	$$\lim\limits_{\nu\to\pm\infty} \underline{\lambda}_k(\nu, \mathcal{A}_1) = \lim\limits_{\nu\to\pm\infty} \underline{\lambda}(\nu, \mathcal{A}_1)=+\infty,
	$$
	which yields that there exist $\nu_k$ and $\nu_0$ such that 
	$$	\underline{\lambda}_k(\nu_k, \mathcal{A}_1)=\inf_{\nu \in \mathbb{R}} \underline{\lambda}_k(\nu, \mathcal{A}_1)\le\inf_{\nu \in \mathbb{R}}  \underline{\lambda}(\nu, \mathcal{A}_1)=\underline{\lambda}(\nu_0, \mathcal{A}_1). $$

	Note that a similar argument as \cite[Step 1 of Proof of Theorem 5.4]{LZh2} yields that \(\underline{\lambda}_k(\nu, \mathcal{A}_1)\nearrow\underline{\lambda}(\nu, \mathcal{A}_1)\)  and is locally Lipschitz continuous with respect to \(\nu\) (one can prove this by similar arguments as \cite[Proposition 2.5]{LZh2}). Therefore, by Dini's theorem, as \(k \to \infty\), 
	$$\min\{\underline{\lambda}_k(\cdot, \mathcal{A}_1),\underline{\lambda}(\nu_0, \mathcal{A}_1)\}\to\min\{\underline{\lambda}(\cdot, \mathcal{A}_1),\underline{\lambda}(\nu_0, \mathcal{A}_1)\}=\underline{\lambda}(\nu_0, \mathcal{A}_1)$$ 
	uniformly on any interval \([a,b]\). 
	In particular,
	$$\lim\limits_{k\to+\infty}\underline{\lambda}_k(\nu_k, \mathcal{A}_1)   \geq\lim\limits_{k\to+\infty}\min\{\underline{\lambda}_k(\nu_k, \mathcal{A}_1),\underline{\lambda}(\nu_0, \mathcal{A}_1)\}=\underline{\lambda}(\nu_0, \mathcal{A}_1).$$
Hence, we deduce that \eqref{eq:star} holds. 

Let us now prove that 	$\inf\limits_{\nu\in\mathbb R}\underline{\lambda}(\nu,\mathcal A_1)\le\lim\limits_{n\to\infty}\overline{\lambda}(0,\mathcal{A}_{\Omega_n}).$
For any given \(\varepsilon > 0\), there exists \(k_0 > 0\) such that for all \(k \ge k_0\),
$$
\inf_{\nu \in \mathbb{R}} \underline{\lambda}_k(\nu, \mathcal{A}_1) \ge \inf_{\nu \in \mathbb{R}} \underline{\lambda}(\nu, \mathcal{A}_1) - \varepsilon/2 .
$$
On the other hand, for any fixed \(k\), according to the proof in the previous part and (ii) of Theorem \ref{thm:main_result2}, we have respectively
$$
\overline{\lambda}_k(0, \mathcal{A}_\mathbb{R}) \ge \inf_{\nu \in \mathbb{R}} \underline{\lambda}_k(\nu, \mathcal{A}_1) ,
$$
and for \(n\) sufficiently large,
$$
\overline{\lambda}_k(0, \mathcal{A}_{\Omega_n}) \ge \overline{\lambda}_k(0, \mathcal{A}_\mathbb{R}) - \varepsilon/2 .
$$
In addition, by the definitions of \(\overline{\lambda}(0, \mathcal{A}_{\Omega_n})\) and \(\overline{\lambda}_k(0, \mathcal{A}_{\Omega_n})\), it follows that
$$
\overline{\lambda}(0, \mathcal{A}_{\Omega_n}) \ge \overline{\lambda}_k(0, \mathcal{A}_{\Omega_n}) \text{ for any }n\in\mathbb Z.
$$

Combining the above inequalities, for any given \(\varepsilon > 0\), there exists \(n_0 > 0\) such that for all \(n \ge n_0\),
\begin{align*}
	&\quad\overline{\lambda}(0, \mathcal{A}_{\Omega_n}) \ge \overline{\lambda}_k(0, \mathcal{A}_{\Omega_n}) \ge \overline{\lambda}_k(0, \mathcal{A}_\mathbb{R}) - \varepsilon/2 \\
	&\ge \inf_{\nu \in \mathbb{R}} \underline{\lambda}_k(\nu, \mathcal{A}_1) - \varepsilon/2 \ge \inf_{\nu \in \mathbb{R}} \underline{\lambda}(\nu, \mathcal{A}_1) - \varepsilon .
\end{align*}
Therefore, by the arbitrariness of \(\varepsilon\), taking the limit inferior on both sides yields:
$$
\liminf_{n \to \infty} \overline{\lambda}(0, \mathcal{A}_{\Omega_n}) \ge \inf_{\nu \in \mathbb{R}} \underline{\lambda}(\nu, \mathcal{A}_1) .
$$
This complete the proof.
\end{proof}

\section{Almost periodic coefficients: proofs of Theorem \ref{thm:main_result3}}

%Similar to the proof of Theorem \ref{thm:main_result1}, we can show that
%\begin{equation}
%	\underline{\lambda}(\nu,\mathcal{A}_1)\leq \underline{\lambda}(\nu,\mathcal{A}_2)\leq \overline{\lambda}(\nu,\mathcal{A}_2)\leq \overline{\lambda}(\nu,\mathcal{A}_1). \label{eq:5.0}
% \end{equation}

To prove Theorem \ref{thm:main_result3}, we need to give and prove the following theorems and lemmas.

\begin{theorem}\label{thm5.1}
   Assume that (J1), (J2), (A1), and (A3) hold, and that \(a(x)\) is uniformly continuous on \(\mathbb{R}\).
   	Let \(\varepsilon>0\) be fixed, and let \(w\in \mathrm{USC}(\mathbb{R})\) and \(v\in \mathrm{LSC}(\mathbb{R})\) be bounded functions. 
   	Suppose that, in the viscosity sense,
	\begin{equation}
		\varepsilon w(x)-\int_{\mathbb{R}}e^{-\nu(x-y)}J(x,y)e^{w(y)-w(x)}\,dy-q(x)w'(x)-a(x)-\nu q(x)\le 0, \label{eq:5.1a}
	\end{equation}
	\begin{equation}
		\varepsilon v(x)-\int_{\mathbb{R}}e^{-\nu(x-y)}J(x,y)e^{v(y)-v(x)}\,dy-q(x)v'(x)-a(x)-\nu q(x)\ge 0. \label{eq:5.1b}
	\end{equation}
	Then \(w\le v\) on \(\mathbb{R}\).
\end{theorem}

\begin{proof}
	
	Define the function
	$$
	\Phi(x,y)=w(x)-v(y)-\alpha|x-y|^2-\mu(x^2+y^2),\quad \alpha>0,\ \mu>0.
	$$
	Clearly
	$$
	w(x)-v(y)-\alpha|x-y|^2-\mu(x^2+y^2)
	\le \sup_{x,y\in\mathbb{R}}\left\{w(x)-v(y)-\alpha|x-y|^2-\mu(x^2+y^2)\right\}.
	$$
	Since the supremum can be attained, there exists \((\xi,\eta)\in\mathbb{R}^2\) such that
	\begin{equation}
		\Phi(x,y)\le \Phi(\xi,\eta)
		=w(\xi)-v(\eta)-\alpha|\xi-\eta|^2-\mu(\xi^2+\eta^2), \label{eq:5.2}
	\end{equation}
	where the point \((\xi,\eta)\) depends on \(\alpha\) and \(\mu\).	
	Taking test functions
	$$
	\phi(x)=v(\eta)+\alpha|x-\eta|^2+\mu(x^2+\eta^2),\text{ and }
	\psi(y)=w(\xi)-\alpha|\xi-y|^2-\mu(\xi^2+y^2).
	$$
	Then \(w(\xi)-\phi(\xi)\) achieves a global maximum, while \(v(\eta)-\psi(\eta)\) attains a global minimum. Set \(\phi_1(x)=\phi(x)+d_1\) with \(d_1=w(\xi)-\phi(\xi)\). Then \(\phi_1(x)\ge w(x)\), \(\phi_1(\xi)=w(\xi)\).
	
	Then we have
	\begin{equation}
		\varepsilon w(\xi)-a(\xi)-\nu q(\xi)-q(\xi)\phi'(\xi)
		-\int_{\mathbb{R}}e^{-\nu(\xi-y)}J(\xi,y)e^{w(y)-w(\xi)}\,dy \le 0. \label{eq:5.3}
	\end{equation}
	
	Similarly, set \(\psi_1(y)=\psi(y)+d_2\) with \(d_2=v(\eta)-\psi(\eta)\), which satisfies \(\psi_1(y)\le v(y)\), \(\psi_1(\eta)=v(\eta)\). Then we have
	\begin{equation}
		\varepsilon v(\eta)-a(\eta)-\nu q(\eta)-q(\eta)\psi'(\eta)
		-\int_{\mathbb{R}}e^{-\nu(\eta-y)}J(\eta,y)e^{v(y)-v(\eta)}\,dy \ge 0. \label{eq:5.4}
	\end{equation}
	
	Suppose that \(w>v\). Then there exist \(x_0\in\mathbb{R}\) and \(\delta>0\) such that
	$$
	w(x_0)-v(x_0)\ge 2\delta.
	$$
	For sufficiently small \(\mu>0\) with \(\mu\le \delta/(2x_0^2)\), we have
	$$
	\delta \le w(x_0)-v(x_0)-\mu(x_0^2+x_0^2)
	\le w(\xi)-v(\eta)-\alpha|\xi-\eta|^2-\mu(\xi^2+\eta^2),
	$$
	which implies \(\delta \le w(\xi)-v(\eta)\). Moreover,
	$$
	\delta+\alpha|\xi-\eta|^2+\mu(\xi^2+\eta^2)
	\le w(\xi)-v(\eta)
	\le \sup_{\mathbb{R}}|w|+\sup_{\mathbb{R}}|v|.
	$$
	Hence, \(\alpha|\xi-\eta|^2\) is bounded and
		\begin{equation}\label{eq:converge}
			\begin{cases}
				\text{\(|\xi-\eta|\to0\) as \(\alpha\to+\infty\) uniformly with respect to \(\mu\),} \\
				\text{\(\mu\xi, \mu\eta\to0\) as \(\mu\to0\) uniformly with respect to \(\alpha\).}
			\end{cases}
		\end{equation}
		Moreover, by the standard property of the penalty method (see \cite[Proposition 3.7]{Crandall1992}), we know that
	$$
	\lim_{\alpha\to+\infty} \alpha|\xi-\eta|^2=0 \text{ uniformly with respect to }\mu.
	$$
	From \eqref{eq:5.3} and \eqref{eq:5.4}, one has
	\begin{equation}
		\begin{aligned}
			\varepsilon\delta
			&\le \varepsilon\bigl(w(\xi)-v(\eta)\bigr) \\
			&\le a(\xi)-a(\eta)+\nu\bigl(q(\xi)-q(\eta)\bigr)
			+q(\xi)\phi'(\xi)-q(\eta)\psi'(\eta) \\
			&\quad +\int_{\mathbb{R}}\Bigl\{
			e^{-\nu(\xi-y)}J(\xi,y)e^{w(y)-w(\xi)}
			-e^{-\nu(\eta-y)}J(\eta,y)e^{v(y)-v(\eta)}
			\Bigr\}\,dy.
		\end{aligned} \label{eq:5.6}
	\end{equation}
	
	To estimate the integral term, we first set
	$$
	C_1:=\max\left\{\exp(\sup\limits_{x}w-\inf\limits_{x}w),\ \exp(\sup\limits_{x}v-\inf\limits_{x}v)\right\}, \text{ and }	C_2:=\int_{\mathbb{R}}J_0(y)e^{-\nu y}\,dy.
	$$
	From assumption (J1), we can choose \(R\) such that
	\begin{equation}
		\int_{|y|>R} J_0(y)e^{-\nu y}\,dy < \frac{\varepsilon\delta}{8C_1}. \label{eq:5.5}
	\end{equation}
	Next, we split the integral in \eqref{eq:5.6} into three parts \(M_1,M_2,M_3\):
	\begin{equation}\label{eq:5.8}
		\begin{aligned}
			&\quad\int_{\mathbb{R}}\Bigl\{e^{-\nu(\xi-y)}J(\xi,y)e^{w(y)-w(\xi)}
			-e^{-\nu(\eta-y)}J(\eta,y)e^{v(y)-v(\eta)}\Bigr\}\,dy \\
			&=\int_{\mathbb{R}}\Bigl\{e^{-\nu s}J(\xi,\xi-s)e^{w(\xi-s)-w(\xi)}
			-e^{-\nu s}J(\eta,\eta-s)e^{v(\eta-s)-v(\eta)}\Bigr\}\,ds \\
			&=\int_{|s|\le R}
			e^{-\nu s}\Bigl[J(\xi,\xi-s)-J(\eta,\eta-s)\Bigr]e^{w(\xi-s)-w(\xi)}\,ds  \\
			&\quad +\int_{|s|\le R}
			\Bigl[e^{w(\xi-s)-w(\xi)}-e^{v(\eta-s)-v(\eta)}\Bigr]J(\eta,\eta-s)e^{-\nu s}\,ds  \\
			&\quad +\int_{|s|>R} e^{-\nu s}J(\xi,\xi-s)e^{w(\xi-s)-w(\xi)}\,ds  -\int_{|s|>R} e^{-\nu s}J(\eta,\eta-s)e^{v(\eta-s)-v(\eta)}\,ds.
		\end{aligned} 
	\end{equation}
	By \eqref{eq:5.5}, we obtain
	\begin{equation}
		\begin{aligned}
			|M_3|
			&\le \int_{|s|>R} e^{-\nu s}J_0(s)e^{w(\xi-s)-w(\xi)}\,ds  +\int_{|s|>R} e^{-\nu s}J_0(s)e^{v(\eta-s)-v(\eta)}\,ds\\
			&\le C_1\frac{\varepsilon\delta}{8C_1}+C_1\frac{\varepsilon\delta}{8C_1}=\frac{\varepsilon\delta}{4}.
		\end{aligned} \label{eq:5.12}
	\end{equation}
	Since assumption (J2) holds and \(|\xi-\eta|\to0\) as \(\alpha\to+\infty\), then
	$$
	e^{-\nu s}\Bigl|J(\xi,\xi-s)-J(\eta,\eta-s)\Bigr|
	< \frac{\varepsilon\delta}{8R\,C_1}.
	$$
	Thus,
	\begin{equation}		 \label{eq:5.9}
		|M_1|\le \int_{|s|\le R}e^{-\nu s}\Bigl|J(\xi,\xi-s)-J(\eta,\eta-s)\Bigr|e^{w(\xi-s)-w(\xi)}\,ds < \frac{\varepsilon\delta}{8}.
		\end{equation}
	Furthermore, from \eqref{eq:5.2}, we get
	$$
	w(\xi-s)-w(\xi) \le v(\eta-s)-v(\eta)+2\mu s^2-2\mu(\xi+\eta)s.
	$$
		Hence,
	\begin{equation}
			M_2\le \int_{|s|\le R}
			\Bigl[e^{2\mu s^2-2\mu(\xi+\eta)s}-1\Bigr]e^{v(\eta-s)-v(\eta)}J(\eta,\eta-s)e^{-\nu s}\,ds \to0\text{ as }\mu\to0. \label{eq:5.11}
	\end{equation}
	
Besides, one can compute the derivative
	$$
	\phi'(\xi)=2\alpha(\xi-\eta)+2\mu\xi,\text{ and }
	\psi'(\eta)=2\alpha(\xi-\eta)-2\mu\eta.
	$$
	Consequently,
	\begin{equation}
		q(\xi)\phi'(\xi)-q(\eta)\psi'(\eta)
		=2\alpha(\xi-\eta)\bigl(q(\xi)-q(\eta)\bigr)
		+2\mu\bigl(\xi q(\xi)+\eta q(\eta)\bigr). \label{eq:5.7}
	\end{equation}
	
	Substituting \eqref{eq:5.12}-\eqref{eq:5.7} into \eqref{eq:5.6}, we deduce that
	\begin{equation}\label{eq:5.14}
		\begin{aligned}
			\varepsilon\delta
			&\le 2\alpha L_q|\xi-\eta|^2
			+\bigl(a(\xi)-a(\eta)\bigr)
			+\nu\bigl(q(\xi)-q(\eta)\bigr) \\
			&\quad +2\mu\bigl(\xi q(\xi)+\eta q(\eta)\bigr)  +M_1+M_2+M_3.
				\end{aligned} 
	\end{equation}
	Recall that \(a(x)\) is uniformly continuous and \eqref{eq:converge}. 
	Taking \(\alpha\to+\infty\) in \eqref{eq:5.14} first and then taking $\mu\to0$ yields
     $$
	\frac{\varepsilon\delta}{2}\le 0,
	$$
	which is a contradiction.	
\end{proof}

\begin{theorem}\label{thm5.2}
	Assume that (J1), (J2), (A1), and (A3) hold, and that \(a(x)\) is uniformly continuous on \(\mathbb{R}\).
	Let \(\varepsilon>0\) be fixed. Then there exists a unique uniformly continuous viscosity solution \(u^\varepsilon\in L^\infty(\mathbb{R})\) of the equation
	\begin{equation}
		\varepsilon u^\varepsilon(x)-\int_{\mathbb{R}}e^{-\nu(x-y)}J(x,y)e^{u^\varepsilon(y)-u^\varepsilon(x)}\,dy
		-q(x)(u^\varepsilon)'(x)-a(x)-\nu q(x)=0 \label{eq:5.15}
	\end{equation}
   with
	$$
	\inf_{x\in\mathbb{R}}\frac{\displaystyle\int_{\mathbb{R}}e^{-\nu(x-y)}J(x,y)\,dy+a(x)+\nu q(x)}{\varepsilon}
	\le u^\varepsilon(x)
	\le \sup_{x\in\mathbb{R}}\frac{\displaystyle\int_{\mathbb{R}}e^{-\nu y}J_0(y)\,dy+a(x)+\nu q(x)}{\varepsilon}.
	$$
    If, furthermore, (A2) is valid, then $u^\varepsilon\in C^1(\mathbb R)$, and  $\sup\limits_{x\in\mathbb R,\varepsilon\in(0,1)}|(u^\varepsilon)^\prime(x)|<+\infty$.
\end{theorem}
	Denote
	$$
	c_1=\inf_{x\in\mathbb{R}}\frac{\displaystyle\int_{\mathbb{R}}e^{-\nu(x-y)}J(x,y)\,dy+a(x)+\nu q(x)}{\varepsilon},\quad
	c_2=\sup_{x\in\mathbb{R}}\frac{\displaystyle\int_{\mathbb{R}}e^{-\nu y}J_0(y)\,dy+a(x)+\nu q(x)}{\varepsilon}.
	$$
	Then $c_1$ and $c_2$ are sub- and super-solutions to \eqref{eq:5.15}.	
	Using Theorem \ref{thm5.1}, we know \(c_1\le u^\varepsilon(x)\le c_2\). 
	Define 
	$$
	A=\Bigl\{u\in \mathrm{USC}(\mathbb{R}):\ c_1\le u(x)\le c_2,\ u \text{ is a viscosity subsolution of } \eqref{eq:5.15}\Bigr\}.
	$$
	In order to prove the existence of a viscosity solution of \eqref{eq:5.15}, we need the following definitions and lemmas.
\begin{definition}\label{def5.3}
	Let \(u:\mathbb{R}\to\mathbb{R}\) be an arbitrary bounded function. The lower (upper) semi-continuous envelope, which we denote by \(u_*\) (\(u^*\)), is given by
	$$
	u_*(x)=\lim_{r\to0^+}\inf_{|y-x|<r}u(y)
	\quad
	\left(
	u^*(x)=\lim_{r\to0^+}\sup_{|y-x|<r}u(y)
	\right).
	$$
\end{definition}
	By definition, there exists a sequence \(x_n\to x_0\) such that \(u(x_n)\to u^*(x_0)\). \(u_*\) is the largest lower semi-continuous function which is smaller than or equal to \(u\). \(u^*\) is the smallest upper semi-continuous function which is larger than or equal to \(u\). Moreover, for any \(x_0\in\mathbb{R}\),\\
	(i) there exists a sequence \(x_n\to x_0\) such that \(\displaystyle\lim_{n\to\infty}u_n(x_n)=u^*(x_0)\).\\
	(ii) for any sequence \(\xi_n\to x_0\), we have
		$
		\limsup\limits_{n\to\infty}u_n(\xi_n)\le u^*(x_0).
		$

\begin{lemma}\label{lem5.4}
  Under the assumptions of Theorem \ref{thm5.2}, the upper semi-continuous envelope \(w_A^*\) of \(w_A=\sup_{u\in A}u(x)\) is a viscosity subsolution of \eqref{eq:5.15}.
\end{lemma}

\begin{proof}
	Let \(\phi\) be a test function that touches \(w_A^*\) at \(x_0\) strictly from above, i.e.,
	$$
	\phi(x_0)=w_A^*(x_0),\quad \phi(x)>w_A^*(x)\quad \text{for }x\in x_0-\mathcal N\setminus\{x_0\},
	$$
	where \(\mathcal N\) is some neighborhood of the origin. First, using Definition \ref{def5.3}, we can find sequences \(x_n\) and \(u_n\in A\) such that
	$$
	\lim_{n\to\infty}x_n=x_0,\qquad \lim_{n\to\infty}u_n(x_n)=w_A^*(x_0).
	$$
	For any sequence \(\xi_n\to x_0\), we have
	\begin{equation}
		\limsup_{n\to\infty}u_n(\xi_n)\le w_A^*(x_0). \label{eq:5.17}
	\end{equation}
	
	Note that \(\phi-w_A^*\) reaches its minimum at \(x_0\) over \(x_0-\mathcal N\). Then one can find that \(\phi-u_n\) reaches a minimum at \(\eta_n\in x_0-\mathcal N/2\) for large \(n\). That is,
	$$
	(\phi-u_n)(\eta_n)\le (\phi-u_n)(x),\quad \text{for all }x\in x_0-\mathcal N/2.
	$$
	Moreover, \(\displaystyle\lim_{n\to\infty}\eta_n=x_0\). In fact, if \(\eta_n\to\bar\eta\in x_0-\mathcal N/2\), then using \eqref{eq:5.17}, we have
	$$
	(\phi-w_A^*)(x_0)=\lim_{n\to\infty}(\phi-u_n)(x_n)
	\ge\liminf_{n\to\infty}(\phi-u_n)(\eta_n)
	\ge\liminf_{\eta\to\bar\eta}(\phi-w_A^*)(\eta)
	\ge(\phi-w_A^*)(\bar\eta).
	$$
	Since \(\phi-w_A^*\) reaches its strict  minimum at \(x_0\) over \(x_0-\mathcal N\), we know 
	$$\eta\to x_0, \text{ and }\displaystyle\lim_{n\to\infty}u_n(\eta_n)=w_A^*(x_0).$$
	
	Recall equation \eqref{eq:5.15}, each \(u_n\) is a viscosity subsolution, so
	\begin{equation}
		\varepsilon u_n(\eta_n)-\int_{\mathbb{R}}e^{-\nu(\eta_n-y)}J(\eta_n,y)\bigl(e^{u_n(y)-u_n(\eta_n)}\bigr)\,dy
		-q(\eta_n)\phi'(\eta_n)-a(\eta_n)-\nu q(\eta_n)\le 0. \label{eq:5.18}
	\end{equation}
	Using Fatou's lemma and the semi-continuity of \(w_A^*\),
	\begin{equation}
		\begin{aligned}
			&\quad\limsup_{n\to\infty}\int_{\mathbb{R}}e^{-\nu(\eta_n-y)}J(\eta_n,y)\bigl(e^{u_n(y)-u_n(\eta_n)}\bigr)\,dy\\
			&\le \int_{\mathbb{R}}\limsup_{n\to\infty}\Bigl(e^{-\nu(\eta_n-y)}J(\eta_n,y)\bigl(e^{u_n(y)-u_n(\eta_n)}\bigr)\Bigr)\,dy \\
			&\le \int_{\mathbb{R}}e^{-\nu(x_0-y)}J(x_0,y)\bigl(e^{w_A^*(y)-w_A^*(x_0)}\bigr)\,dy.
		\end{aligned} \label{eq:5.19}
	\end{equation}
	With \eqref{eq:5.19}, take \(n\to\infty\) in \eqref{eq:5.18} to deduce that
	$$
	\varepsilon w_A^*(x_0)-\int_{\mathbb{R}}e^{-\nu(x_0-y)}J(x_0,y)\bigl(e^{w_A^*(y)-w_A^*(x_0)}\bigr)\,dy
	-q(x_0)\phi'(x_0)-a(x_0)-\nu q(x_0)\le 0.
	$$
	Since \(\phi\) touches \(w_A^*\) from above at \(x_0\) over \(x_0-\mathcal N\), we conclude that \(w_A^*\) is a viscosity subsolution to \eqref{eq:5.15}.
\end{proof}

\begin{lemma}\label{lem5.5}
	Assume that (J1), (J2), and (A1) hold.
	Let \(u\) be a bounded viscosity subsolution of \eqref{eq:5.15}. Suppose that \(u_*\) is not a viscosity supersolution. Then there exists a test function \(\psi\) touching \(u_*\) from below at some point \(x_0\) over \(x_0-\mathcal N\), such that
	\begin{equation}
		\varepsilon u_*(x_0)-\int_{\mathbb{R}}\bigl(e^{u_*(y)-u_*(x_0)}\bigr)e^{-\nu(x_0-y)}J(x_0,y)\,dy
		-q(x_0)\psi'(x_0)-a(x_0)-\nu q(x_0)\le -\eta, \label{eq:5.20}
	\end{equation}
	for some \(\eta>0\). Then there exists an upper semicontinuous function \(\widetilde u\in \mathrm{USC}(\mathbb{R})\) such that \(\widetilde u\ge u\), \(\widetilde u\neq u\), and \(\widetilde u\) is a viscosity subsolution.
\end{lemma}

\begin{proof}
	Without loss of generality, we may assume \(x_0=0\) and that a strict minimum is attained at this point. 
	Denote $$E(x):= \varepsilon\psi(x) - \int_{\mathbb{R}} (e^{u_*(y)-\psi(x)})e^{-\nu(x-y)}J(x,y)dy - q(x)\psi'(x) - a(x) - \nu q(x).$$
	Then \eqref{eq:5.20} holds with \(x_0=0\),
	\begin{equation}
		E(0)=\varepsilon u_*(0)-\int_{\mathbb{R}}\bigl(e^{u_*(y)-u_*(0)}\bigr)e^{-\nu(0-y)}J(0,y)\,dy
		-q(0)\psi'(0)-a(0)-\nu q(0)\le -\eta. \label{eq:5.21}
	\end{equation}
	By taking \(r_0>0\) sufficiently small, we can ensure that, for any $ x\in B_{r_0}$,
	\begin{equation}\label{eq:5.22}
		|E(x)-E(0)|\leq\frac{\eta}{4},
	\end{equation}
	since  (J1), (J2), and (A1) hold.
	Define
	\[
	U(x)=
	\begin{cases}
		\psi+\delta_r, & \text{ on } B_{r_0},\\4pt]
		u_*(x), & \text{ on } \mathbb{R}\setminus (-B_{r_0}),
	\end{cases}
	\]
	where \(\delta_r>0\) depending on $r<r_0$ will be chosen later. For sufficiently small $r_0$ and \(\delta_r\), we have	
	\begin{equation}\label{eq:5.23}
		\begin{cases}
			\left|\int_{B_{r_0}}\bigl(e^{u_*(y)-\psi(x)}-e^{\psi(y)-\psi(x)}\bigr)e^{-\nu(x-y)}J(x,y)\,dy\right|<\frac{\eta}{4},\\
			\varepsilon\delta_r+\left|\int_{B_{r_0}^c}\bigl(e^{u_*(y)-\psi(x)}-e^{u_*(y)-\psi(x)-\delta_r}\bigr)e^{-\nu(x-y)}J(x,y)\,dy\right|<\frac{\eta}{4}. 
		\end{cases}
			\end{equation}
	Note that \(\psi(0)=u_*(0)\). Taking \(r_0>0\) small enough, and  combining \eqref{eq:5.21}, \eqref{eq:5.22}, and \eqref{eq:5.23}, we obtain that, for all \(x\in B_{r_0}\):
\begin{equation}\label{eq:5.24} 
	\begin{aligned}
		 &\quad \varepsilon U(x) - \int_{\mathbb{R}} (e^{U(y)-U(x)})e^{-\nu(x-y)}J(x,y)dy - q(x)\psi'(x) - a(x) - \nu q(x) \\
		&= \varepsilon(\psi(x) + \delta_r) - \int_{\mathbb{R}} (e^{U(y)-\psi(x)-\delta_r})e^{-\nu(x-y)}J(x,y)dy - q(x)\psi'(x) - a(x) - \nu q(x) \\
		&= E(x) + \varepsilon\delta_r + \int_{\mathbb{R}} (e^{u_*(y)-\psi(x)} - e^{U(y)-\psi(x)-\delta_r})e^{-\nu(x-y)}J(x,y)dy \\
		&= E(x) +  \varepsilon\delta_r + \int_{B_{r_0}^c} \left( e^{u_*(y)-\psi(x)} - e^{u_*(y)-\psi(x)-\delta_r} \right) e^{-\nu(x-y)}J(x,y)dy  \\
		&\quad + \int_{B_{r_0}} \left( e^{u_*(y)-\psi(x)} - e^{\psi(y)-\psi(x)} \right) e^{-\nu(x-y)}J(x,y)dy \\
		&< -\frac{3\eta}{4} + \frac{\eta}{4} + \frac{\eta}{4} = -\frac{\eta}{4}.
	\end{aligned}
\end{equation}
		Choose \(\delta_r<\min\{u_*(x)-\psi(x):x\in\overline{B_{r_0}}\setminus\overline{B_{r/2}}\}\). Then \(U(x)<u_*(x)\le u(x)\) for \(x\in\overline{B_{r_0}}\setminus\overline{B_{r/2}}\), which yields that the function \(\widetilde u=\max\{u,U\}\) belongs to $\mathrm{USC}(\mathbb{R})$. We claim that \(\widetilde u\) is a viscosity subsolution of \eqref{eq:5.15}.
	
	Suppose \(\phi\) is a test function touching \(\widetilde u\) from above at some point \(\xi\), so that \(\phi(\xi)=\widetilde u(\xi)\) and \(\phi(x)>\widetilde u(x)\) for $x\in$ some neighborhood $\mathcal N$ of $\xi$. There are two cases:
	
	\noindent\textbf{Case 1:} \(\widetilde u(\xi)=u(\xi)\). Then
	$$
	\begin{aligned}
		&\quad\varepsilon\widetilde u(\xi)-\int_{\mathbb{R}}\bigl(e^{\widetilde u(y)-\widetilde u(\xi)}\bigr)e^{-\nu(\xi-y)}J(\xi,y)\,dy
		-q(\xi)\phi'(\xi)-a(\xi)-\nu q(\xi) \\
		&\le\varepsilon u(\xi)-\int_{\mathbb{R}}\bigl(e^{u(y)-u(\xi)}\bigr)e^{-\nu(\xi-y)}J(\xi,y)\,dy
		-q(\xi)\phi'(\xi)-a(\xi)-\nu q(\xi) \le 0,
	\end{aligned}
	$$
	since \(u\) is a viscosity subsolution.
	
	\noindent\textbf{Case 2:} \(\widetilde u(\xi)=U(\xi)\). Note that if \(x\notin B_{r/2}\), then \(\widetilde u(x)=u(x)\). We therefore assume that \(\xi\in B_{r/2}\). 
		Since \(\phi(\xi)=\widetilde u(\xi)= U(\xi)=\psi(\xi)+\delta_r\), and \(\phi(x)\ge\widetilde u(\xi)\ge U(x)\ge \psi(x)+\delta_r\) for \(x\in B_{r_0}\cap\mathcal N\), we have \(\phi'(\xi)= \psi'(\xi)\). 
	Then, using \eqref{eq:5.24},
	$$
	\begin{aligned}
		&\quad\varepsilon\widetilde u(\xi)-\int_{\mathbb{R}}\bigl(e^{\widetilde u(y)-\widetilde u(\xi)}\bigr)e^{-\nu(\xi-y)}J(\xi,y)\,dy
		-q(\xi)\phi'(\xi)-a(\xi)-\nu q(\xi) \\
		&\le\varepsilon U(\xi)-\int_{\mathbb{R}}\bigl(e^{U(y)-U(\xi)}\bigr)e^{-\nu(\xi-y)}J(\xi,y)\,dy
		-q(\xi)\phi'(\xi)-a(\xi)-\nu q(\xi)\\
		&= \varepsilon U(\xi)-\int_{\mathbb{R}}\bigl(e^{U(y)-U(\xi)}\bigr)e^{-\nu(\xi-y)}J(\xi,y)\,dy
		-q(\xi)\psi'(\xi)-a(\xi)-\nu q(\xi) \le -\frac{\eta}{4}.
	\end{aligned}
	$$

Combining both cases, \(\widetilde u\) is a viscosity subsolution.	
	By construction, \(\widetilde u\ge u\). It remains to show that \(\widetilde u\neq u\). Take a sequence \(x_n\to 0\) such that \(\displaystyle\lim_{n\to\infty}u(x_n)=u_*(0)\). Then
	$$
	\lim_{n\to\infty}\bigl(\widetilde u(x_n)-u(x_n)\bigr)
	=\lim_{n\to\infty}\bigl(\psi(x_n)+\delta_r-u(x_n)\bigr)
	=\psi(0)+\delta_r-u_*(0)=\delta_r>0.
	$$
	Hence for sufficiently large \(n\), \(\widetilde u(x_n)>u(x_n)\), so \(\widetilde u\neq u\).
\end{proof}
 Next, we prove Theorem \ref{thm5.2}.
\begin{proof}[Proof of Theorem \ref{thm5.2}]
	We first show that equation \eqref{eq:5.15} admits a unique bounded viscosity solution \(u^\varepsilon\).
	By Lemma \ref{lem5.4}, we see that \(w_A^*=w\) is a viscosity subsolution of \eqref{eq:5.15}. We claim that \(w_*\) is a viscosity supersolution of \eqref{eq:5.15}. 
	If not, by Lemma \ref{lem5.5}, we can find a viscosity subsolution \(\widetilde w\ge w\) such that \(\widetilde w(x)>w(x)\) at some points, which contradicts the definition of \(w_A^*\). 
	Therefore, \(w_*\) is a viscosity supersolution of \eqref{eq:5.15}. Moreover, by Theorem \ref{thm5.1}, we have \(w_A^*\le w_*\). On the other hand, \(w_A^*\ge (w_A^*)_*=w_*\). Hence, \(w_A^*=w_*\) is the unique bounded viscosity solution to \eqref{eq:5.15}. We denote this solution by $u^\varepsilon$.
	
	Let us now show that  \(u^\varepsilon\) is a \(C^1\) function, and $\sup\limits_{x\in\mathbb R,\varepsilon\in(0,1)}|(u^\varepsilon)^\prime(x)|<+\infty$ with (A3).
		Define
	\begin{equation}
		g(x)=\frac{\varepsilon u^\varepsilon(x)-a(x)-\nu q(x)-\displaystyle\int_{\mathbb{R}}e^{-\nu(x-y)}J(x,y)\bigl(e^{u^\varepsilon(y)-u^\varepsilon(x)}\bigr)\,dy}{q(x)}. \label{eq:5.25}
	\end{equation}
	Then $g$ is bounded continuous function, and there exists a differentiable function \(G(x)\) on \(\mathbb{R}\) such that
	$$
	G'(x)=g(x).
	$$
	Let 
		\begin{equation*}
		F:=\begin{cases}
			u^\varepsilon-G  \text{ if }q>0,\\
			G-u^\varepsilon  \text{ if }q<0. 
		\end{cases}
	\end{equation*}
	Then one can check that for any test function \(\phi\in C^1\) touch $F$ from above (resp. below) at $x_0$, then $\phi^\prime(x_0)\ge0$ (resp. $\le0$).
    Consider the sup- and inf-convolution approximation $F^\theta$ and $F_\theta$  of $F$.
    We can prove the following in a similar manner to the claim in the proof of Theorem \ref{thm:main_result1}:
    \begin{equation*}
    	\begin{cases}
    		F^{\theta,\prime}(x)\ge0, \text{ a.e. }x\in\mathbb R,\\
    		F_\theta^\prime(x)\leq0, \text{ a.e. }x\in\mathbb R. 
    	\end{cases}
    \end{equation*}
  Since both \(F_\theta\) and \(F^\theta(x)\) converge uniformly to \(F\) by (i) of Proposition \ref{prop:convolution}, we deduce that $F$ is a constant function.
  Therefore, $u^\varepsilon$ is differentiable and is the unique bounded classical solution to \eqref{eq:5.15}.

	To show $\sup\limits_{x\in\mathbb R,\varepsilon\in(0,1)}|(u^\varepsilon)^\prime(x)|<+\infty$, we only need to prove
	$$
	\sup\limits_{x\in\mathbb R,\varepsilon\in(0,1)}\int_{\mathbb{R}}e^{-\nu(x-y)}J(x,y)e^{u^\varepsilon(y)-u^\varepsilon(x)}\,dy<+\infty.
	$$
Note that equation \eqref{eq:5.15} can be transformed into
	$$
	\int_{\mathbb{R}}e^{-\nu(x-y)}J(x,y)e^{u^\varepsilon(y)}\,dy+\bigl(a(x)+\nu q(x)-\varepsilon u^\varepsilon(x)\bigr)e^{u^\varepsilon(x)}+q(x)\left(e^{u^\varepsilon(x)}\right)^\prime=0.
	$$
	 Recall that $|\varepsilon u^\varepsilon|_\infty$ is uniformly bounded with respect to $\varepsilon$). Then we apply  Theorem \ref{harnack} to $e^{u^\varepsilon(x)}$ that there exists $C>0$ such that (without loss of generality, assume $\xi<\eta$)
	 $$
	 e^{u^\varepsilon(\xi)-u^\varepsilon(\eta)}=\frac{e^{u^\varepsilon(\xi)}}{e^{u^\varepsilon(\xi+1)}}
	 \cdot \frac{e^{u^\varepsilon(\xi+1)}}{e^{u^\varepsilon(\xi+2)}}
	 \cdots
	 \frac{e^{u^\varepsilon(\xi+[\eta-\xi])}}{e^{u^\varepsilon(\eta)}}
	 \le C^{\left( \left| \eta-\xi\right| + 1 \right)}.
	 $$
	 Consequently,
	 \begin{equation}\label{eq:linear growth1}
	 	|u^{\varepsilon}(x)-u^{\varepsilon}(y)| 	\le C\left( \left| {y-x} \right| + 1 \right).
	 \end{equation}
	which yields that 
		$$
	\begin{aligned}
		&\quad\sup\limits_{x\in\mathbb R,\varepsilon\in(0,1)}\int_{\mathbb{R}}e^{-\nu(x-y)}J(x,y)e^{u^\varepsilon(y)-u^\varepsilon(x)}\,dy \\
		&\le\sup\limits_{x\in\mathbb R,\varepsilon\in(0,1)}\int_{\mathbb{R}}e^{-\nu(x-y)}J(x,y)e^{C\left( \left| {y-x} \right| + 1 \right)}\,dy
		\le e^C\int_{\mathbb{R}}e^{(C-\nu)y}J_0(y)\,dy.
	\end{aligned}
	$$
	Hence $(u^\varepsilon)^\prime=g$ is bounded.
\end{proof}

\begin{lemma}\label{lem5.9}
	Assume that (J1), (J2)$^\prime$,  (A1)$^\prime$, (A2), and (A3) hold. Then the limit \(\lim\limits_{\varepsilon\to0}\varepsilon u^\varepsilon\)
	exists in $C(\mathbb R)$, where \(u^\varepsilon\in L^\infty(\mathbb{R})\) is the viscosity solution to equation \eqref{eq:5.15}.
\end{lemma}

\begin{proof}
For any \(x_n, x_m \in \R\), we set the following notation:
$$
u^\varepsilon_n(x) = u^\varepsilon(x + x_n), \ a_n(x) = a(x+x_n), \ q_n(x) = q(x+x_n), \ J_n(x,z) = J(x+x_n, x+x_n-z),
$$
and similarly for the index \(m\).
Let
$$
v(x):=u^\varepsilon_n(x) + \frac{T(x_n, x_m)}{\varepsilon},
$$
where
\begin{align*}
	T(x_n, x_m) &:= \|a_n - a_m\|_\infty + \|(u^\varepsilon)'\|_\infty \|q_n - q_m\|_\infty \\
	&\quad + \sup_x \int_{\R} |J_n(x,z) - J_m(x,z)| e^{(L_q+\nu)|z|} \,dz.
\end{align*}
From the equations satisfied by \(u^\varepsilon_n\) and \(u^\varepsilon_m\), we have
\begin{equation*}
	\varepsilon u^\varepsilon_m(x) - a_m(x) - \nu q_m(x) - \int_{\R} e^{-\nu z} J_m(x,z) e^{u^\varepsilon_m(x-z) - u^\varepsilon_m(x)} \,dz - q_m(x)(u^\varepsilon_m)'(x) = 0,
\end{equation*}
and
\begin{align*}
	&\quad\varepsilon v(x) - a_m(x) - \nu q_m(x) - \int_{\R} e^{-\nu z} J_m(x,z) e^{u^\varepsilon_n(x-z) - u^\varepsilon_n(x)} \,dz - q_m(x)(u^\varepsilon_n)'(x) \\
	&= T(x_n, x_m) + a_n(x) - a_m(x) + \nu \bigl(q_n(x) - q_m(x)\bigr) \\
	&\quad + \int_{\R} e^{-\nu z} \bigl( J_n(x,z) - J_m(x,z) \bigr) e^{u^\varepsilon_n(x-z) - u^\varepsilon_n(x)} \,dz \\
	&\ge T(x_n, x_m) - \|a_n - a_m\|_\infty - \|(u^\varepsilon)'\|_\infty \|q_n - q_m\|_\infty \\
	&\quad - \int_{\R} |J_n(x,z) - J_m(x,z)| e^{L|z|} e^{|\nu||z|} \,dz \\
	&\ge 0.
\end{align*}
By Theorem \ref{thm5.1}, we obtain
$$
u^\varepsilon_m(x) \le u^\varepsilon_n(x) + \frac{T(x_n, x_m)}{\varepsilon}.
$$
Hence for any \(x \in \R\), by symmetry, we get
\begin{equation}
|\varepsilon u^\varepsilon(x+x_n) - \varepsilon u^\varepsilon(x+x_m)|=	|\varepsilon u^\varepsilon_n(x) - \varepsilon u^\varepsilon_m(x)| \le T(x_n, x_m). \tag{5.22}
\end{equation}
Let
$$
\hat{u}^\varepsilon(x): = u^\varepsilon(x) - u^\varepsilon(0).
$$
In particular, taking \(x=0\) yields
\begin{equation}
	|\varepsilon \hat{u}^\varepsilon(x_n) - \varepsilon \hat{u}^\varepsilon(x_m)| \le T(x_n, x_m). \tag{5.23}
\end{equation}

\vspace{1em}
\noindent\textbf{Claim:} \(\varepsilon \hat{u}^\varepsilon(x) \to 0\) as \(\varepsilon \to 0\) uniformly with respect to \(x\).

\vspace{0.5em}
\noindent\textit{Proof of the claim:} Assume by contradiction that there exist \(\varepsilon_n \to 0\), \(x_n\), \(\theta > 0\) such that
$$
|\varepsilon_n \hat{u}^{\varepsilon_n}(x_n)| \ge 2\theta.
$$
We will prove the case where \(\varepsilon_n \hat{u}^{\varepsilon_n}(x_n) \ge 2\theta\), and one can prove the case where \(\varepsilon_n \hat{u}^{\varepsilon_n}(x_n) \le -2\theta\) similarly.
Note that (J2)$^\prime$ and (A2)$^\prime$ hold.
Therefore, by passing to a subsequence (still denoted by \(\{x_n\}\)), we have
$$
\|a_n - a_m\|_\infty + \|(u^\varepsilon)'\|_\infty \|q_n - q_m\|_\infty \to 0 \quad \text{as } m, n \to \infty.
$$
Furthermore, for any \(\theta > 0\), there exists \(R > 0\) such that
$$
\int_{|z| > R} J_0(z) e^{\alpha|z|} \,dz < \frac{\theta}{8}.
$$
For this fixed \(R\), by (J2)$^\prime$, there exists \(N \in \mathbb{N}\) such that for all \(m, n > N\),
$$
e^{\alpha R} \sup_x \int_{|z| \le R} |J_n(x,z) - J_m(x,z)| \,dz < \frac{\theta}{2}.
$$
Consequently, for \(m, n > N\) (enlarge $N$ if necessary), we obtain
\begin{align*}
	T(x_n, x_m) &\le \|a_n - a_m\|_\infty + \|(u^\varepsilon)'\|_\infty \|q_n - q_m\|_\infty \\
	&\quad + \sup_x \int_{\R} |J_n(x,z) - J_m(x,z)| e^{\alpha|z|} \,dz \\
	&\le \frac{\theta}{4} + e^{\alpha R} \sup_x \int_{|z| \le R} |J_n(x,z) - J_m(x,z)| \,dz + 2 \int_{|z| > R} J_0(z) e^{\alpha|z|} \,dz \\
	&< \frac{\theta}{4} + \frac{\theta}{2} + 2 \cdot \frac{\theta}{8}=\theta.
\end{align*}
We choose and fix a sufficiently large integer \(m_0 \ge N\). Then for all \(n\) sufficiently large, 
\begin{align*}
	2\theta &\le \varepsilon_n \hat{u}^{\varepsilon_n}(x_n) \le \varepsilon_n \hat{u}^{\varepsilon_n}(x_{m_0}) + |\varepsilon_n \hat{u}^{\varepsilon_n}(x_n) - \varepsilon_n \hat{u}^{\varepsilon_n}(x_{m_0})| \\
	&\le \varepsilon_n \hat{u}^{\varepsilon_n}(x_{m_0}) + T(x_n, x_{m_0})< \varepsilon_n \|(u^{\varepsilon_n})^\prime\|_{\infty}|x_{m_0}| + \theta 
\end{align*}
Letting \(n\to \infty\), since \(\varepsilon_n \to 0\) and \(x_{m_0}\) is fixed, we obtain
$$
2\theta \le \theta,
$$
which is a contradiction since \(\theta > 0\). Then, the proof of claim is complete.

The claim means that for any sequence \(\{\varepsilon_n\}\) there exists a subsequence still denoted by \(\{\varepsilon_n\}\) such that \(\varepsilon_n u^{\varepsilon_n} \to \lim_{n\to\infty} \varepsilon_n u^{\varepsilon_n}(0)\) in \(C(\R)\). Then, follow the same idea as in the proof of Step 5 in \cite[Theorem 3.10]{LZh2}, we can carry out that for any sequence \(\varepsilon_n\) tending to \(0\), \(\varepsilon_n u^{\varepsilon_n}\) converges to the same constant as \(n \to +\infty\).
\end{proof}
 we denote the uniformly limit of $\varepsilon u^\varepsilon$ by $\lambda(\nu)$.
\begin{proof}[Proof of Theorem \ref{thm:main_result3}]
	Theorem \ref{thm5.2} proves that the viscosity solution \(u^\varepsilon\) of \eqref{eq:5.15} exists and is a classical solution, and Lemma \ref{lem5.9} proves that \(\varepsilon u^\varepsilon\) converges uniformly to \(\lambda\). Therefore, for any \(\eta > 0\), we have 
	$$
	(\lambda-\eta)e^{u^\varepsilon}\le\mathcal L e^{u^\varepsilon}=\varepsilon u^\varepsilon e^{u^\varepsilon}\le(\lambda+\eta)e^{u^\varepsilon}.
	$$
	 Taking \(e^{u^\varepsilon}\) as a test function and using the arbitrariness of \(\eta\), we immediately obtain 
	 \begin{equation}\label{}
	 	 \underline{\lambda}(\nu,\mathcal{A}_1) = \underline{\lambda}(\nu,\mathcal{A}_2) = \overline{\lambda}(\nu,\mathcal{A}_2) = \overline{\lambda}(\nu,\mathcal{A}_1)=\lambda(\nu).
	 \end{equation}
	The last assertion of Theorem 1.4 follows by combining (iv) of Theorem \ref{thm:main_result2} with above equation.
		\end{proof}

%\section{Appendix}
\section*{Conflict of interest}
The authors declare that there is no conflict of interest.

\section*{Use of AI tools}

The human authors formulated the mathematical problem, conceived and developed the main results, and independently verified all conclusions. 
AI assisted in improving the language and presentation. The final manuscript was approved by the authors, who take full responsibility for its content.

\end{document}